\documentclass[12pt]{amsart}
\newtheorem{theorem}{Theorem}
\newtheorem{lemma}[theorem]{Lemma}
\newtheorem{proposition}[theorem]{Proposition}
\newtheorem{corollary}[theorem]{Corollary}
\newtheorem{definition}[theorem]{Definition}
\newtheorem{remark}[theorem]{Remark}
\newtheorem{claim}[theorem]{Claim}
\newtheorem{condition}[theorem]{Condition}

\newcommand{\bR}{ {\bf R}}
\newcommand{\dist}{ \mbox{dist}}
\newcommand{\divg}{ \mbox{div}}
\newcommand{\diam}{ \mbox{diam}}
\newcommand{\Lip}{ \mbox{Lip}}

\begin{document}

\title{Variational evolution problems and nonlocal geometric motion}

\author{Mikhail Feldman}
\address{\hskip-\parindent
        Mathematical Sciences Research Institute\\
        1000 Centennial Dr.\\
        Berkeley CA 94720}
\email{mikhail@msri.org}

\let\cal\mathcal

\thanks{Research is supported in part by NSF grants DMS-9701755 (MSRI)
and DMS-9623276.}


\begin{abstract}
We consider two variational evolution problems related to 
Monge-Kantorovich mass transfer. These problems provide
models for collapsing sandpiles and for compression molding.
We prove the following connection between these problems and 
nonlocal geometric curvature motion: The distance functions
to surfaces moving according to certain nonlocal geometric 
laws are solutions of the  variational evolution problems.
Thus we do the first step of the proof of heuristics 
developed in earlier works. The main techniques we use
are differential equations methods in the Monge-Kantorovich theory.
\end{abstract}

\maketitle

In this paper we study two models
involving limits as $p \rightarrow \infty$ of solutions
of $p$-Laplacian evolution problems.
One is a model of collapsing sandpiles proposed by 
L. C. Evans, R. F. Gariepy and the author \cite{EFG}.
 Another is a
model of compression molding proposed by G. Aronsson and 
L. C. Evans  \cite{AEComprMot}.
In both models the limits were characterized as solutions
 of variational evolution problems related to Monge-Kantorovich
mass transfer.
Solutions that
have the form of the distance function to a moving boundary
arise naturally in the both models.
The equations of the motion of the boundary
for both models
were derived {\it heuristically} in \cite{EFG} and \cite{AEComprMot}.
According to the equations
the outer normal velocity of the boundary at a point depends
on both the local geometry (curvatures) and the nonlocal geometry
of the boundary. Thus, the motion
of the boundary is a nonlocal geometric motion.

In this paper we prove rigorously the connection
between geometric and variational evolution problems. Namely,
assuming that a moving surface satisfying the 
geometric equation is given, we prove that the distance
function defines a solution of the corresponding variational
evolution problem. The main assumption is that the surface
remains convex (or, more generally, semiconvex) during the
evolution. For the collapsing sandpiles model we also prove 
the corresponding result for the solutions that have
the form of the maximum of several distance functions
(such solutions represent several sand cones interacting in the
process of collapse). In the proofs we utilize the connection
between the models and the 
Monge-Kantorovich mass transfer. 
This allows us to use the differential equations methods in
the Monge-Kantorovich theory, which have been developed recently 
by L. C. Evans and W. Gangbo \cite{EG}.

The examples suggest that solutions of the geometric
equations derived in  \cite{EFG} can 
develop singularities, even if the initial data is smooth.
Thus we do not assume below
 that the moving surface is smooth. This however makes 
the technique more involved.

In the forthcoming paper \cite{Feld} we construct convexity preserving
viscosity solutions of the geometric equations given convex 
initial
data, thus finishing the proof of the heuristics
developed in \cite{EFG} and \cite{AEComprMot}.

This paper is organized as following.
In Sections \ref{ColSand} - \ref{Th.4} we work on the
collapsing sandpiles model.
In Section \ref{CompMod} we work on the compression molding model.
In Appendix \ref{Apndx 1} we estimate the local
Lipschitz constant of the gradient of the distance function 
of a set.

\section{Notation}
\label{NotatSect}

Let $\Omega \in \bR^n$ be a bounded open set. Let $\Gamma = \partial \Omega$.
We will use the following two distance functions of the set $\Omega$.
The first, interior distance function, is 
\begin{equation}
\label{distFunc00}
d_{\Omega}(x)=
\left\{  
\begin{array}{ll}
\dist(x, \partial \Omega),  & x \in \Omega, \\
0, & x\in \bR^n \setminus \Omega.
\end{array}
\right.
\end{equation}
The second, signed distance function, is
\begin{equation}
\label{sgnDist}
d^{s}_{\Omega} (x)=\left\{
\begin{array}{cr}\mbox{dist}(x,\partial \Omega), & x\in \Omega,
\\ -{\mbox{dist}}(x,\partial \Omega), & x\notin \Omega.
\end{array}
\right. 
\end{equation}

Let $x \in R^n$. 
We denote ${\cal N}_{
\partial \Omega}(x)$
the set of all points of $\partial \Omega$ nearest to $x$,
i.e.,
$$
{\cal N}_{\partial \Omega}(x) = \{
y \in \partial \Omega, \;\;  \;\;
|x-y| = \inf_{z \in \partial \Omega} |x-z| \}.
$$
Note that ${\cal N}_{\partial \Omega}(x)$ may have more
then one point.

Let $y \in \partial \Omega$. Then $\partial \Omega$ is 
differentiable at $y$ if
there exists $\delta>0$ such that
 in a suitable coordinate system 
$(x_1, ..., x_n)$ in $\bR^n$ we have
\begin{equation}
\label{SurfDiff}
 y=0,\; \;\; \Omega
\cap
B_\delta(0)=\{ x_n>\Psi(x')\} \cap
B_\delta(0) \;\;\; \mbox{ where } \;\; \Psi(x')=
o(|x'|).
\end{equation}
Here 
$x'=(x_1, ..., x_{n-1})$, and 
$B_\delta(0) = B^n_\delta(0)$ is the ball in $\bR^n$ 
with 
radius $\delta$ and center at $0$.
$\partial \Omega$ is twice 
differentiable at $y$ if in a 
suitable coordinate system in $\bR^{n-1}$
the function $\Psi$ has the expression 
\begin{equation}
\label{SurfTwiceDiff}
\Psi(x')=
\sum^{n-1}_{i=1} \kappa_ix^2_i+o(|x'|^2).
\end{equation}
The numbers $\kappa_1,...,$
$\kappa_{n-1}$ are the principal curvatures of $\partial \Omega$ at
$y$.

 Let $\{\Omega_t \},\; t\in \bR^1_+$, be a 
family of open bounded sets in $\bR^n$.
 We call a family of sets  $\{\Omega_t \}$ Lipschitz (with respect to
$t$) if
there exists a constant $M>0$ such that for any $t_1, t_2$ 
the set 
$(\Omega_{t_1} \setminus\Omega_{t_2}) \cup 
(\Omega_{t_2} \setminus\Omega_{t_1}) $
is subset of both
$M\mid t_1 - t_2 \mid$ neighborhood (in $\bR^n$) of
 $\partial \Omega_{t_1}$ and $M\mid t_1 - t_2 \mid$ neighborhood of
$\partial \Omega_{t_2}$.

Let 
\begin{eqnarray}
&& E=\cup_t \left(\Omega_t \times \{ t \}\right)
\subset \bR^n \times \bR^1_+, 
\; \;\;\;
\nonumber \\
&& {\bf \Gamma} = \cup_t \left(\partial \Omega_t \times \{ t \}\right)
\subset \bR^n \times \bR^1_+, \;\;\;
\label{defE} \\
&& \Gamma_t = \partial \Omega_t \subset \bR^n.
\nonumber
\end{eqnarray}

Let $(y,t) \in \bR^n \times \bR^1_+$ be such that
$y \in \Gamma_t$ and the function
$(x,\tau) \rightarrow d^{s}_{\Omega_{\tau}} (x)$ is
differentiable at $(y,t)$. We define the outer normal
velocity of  $\Gamma _t$ at $y$ as 
\begin{equation}
\label{outNormVelDef}
V(y) = {\partial \over \partial
t}d^{s}_{\Omega_t} (y).
\end{equation}

Let $(y,t_0) \in {\bf \Gamma}$. The surface ${\bf \Gamma}$  
is (2,1) differentiable
at $(y,t_0)$ if in a suitable coordinate system 
$(x_1, ..., x_n)$ in $\bR^n$ we have
$ y=0,\; \; E
\cap
B^{n+1}_\delta(0)=\{ x_n>\Psi(x', t)\} \cap
B^{n+1}_\delta(0)$ for small $\delta>0$ and 
$$
 \Psi(x',t)= v(t-t_0) + \sum^{n-1}_{i=1} \kappa_ix^2_i 
+ o(\mid x'-y'\mid^2 + \mid t-t_0 \mid) .
$$

${\cal L}^n$ is $n$-dimensional Lebesgue measure.

${\cal H}^n$ is $n$-dimensional Hausdorff measure.

\section{The collapsing sandpiles model}
\label{ColSand}

In the paper \cite{EFG} the  
$p$-Laplacian evolution problems
\begin{equation}\label{pHeatUnS}
\left\{  
\begin{array}{ll}
\partial_t u_{p} - \mbox{div}(\mid Du_p \mid ^ {p-2} Du_p ) =0, & 
\mbox{in}\;\; \bR^n \times (0, \infty), \\
u_p=g, & \mbox{on}\;\; \bR^n \times \{0\},
\end{array}
\right.
\end{equation}
are considered in the ``infinitely fast/slow diffusion
limit'' $p \rightarrow \infty$. 
$Du_p$ denotes
the gradient of $u_p$ with respect to the spatial variables
$x_1, ..., x_n$. The operator
$$
\Delta_p(u) = \mbox{div}(\mid Du_p \mid ^ {p-2} Du_p )
$$
is $p$-Laplacian.
The
initial data $g$ is nonnegative,
Lipschitz continuous function with
compact support. We assume that 
\begin{equation}
\label{UnstLip}
\Lip[g] = L > 1
\end{equation}
where $\Lip[g]$ is the Lipschitz constant of $g$.

 The initial data $g$ satisfying (\ref{UnstLip})
is unstable in the following sense.
It is shown in \cite{EFG} that $u_p \rightarrow u$ uniformly on 
$\bR^n \times (0, \infty)$, and 
that 
\begin{equation}
\label{UnstLipU}
\mid Du \mid \leq 1 \;\;\; \mbox{ a.e. in } \;\; \bR^n \times (0, \infty).
\end{equation}
This limit $u$ does not
depend on $t$ for $t>0$, i.e., $u=u(x)$. By (\ref{UnstLip})
and (\ref{UnstLipU}) we get $u \neq g$. Thus the limit
solution is discontinuous at $t=0$. The 
 transformation $g \rightarrow u$ takes place at $t=0$.
 
The described  problem can be interpreted as a crude model
of a collapsing sandpile. The function $u(x,t)$ is the 
height of the pile at the location $x$, time $t$. 
The main physical assumption is that a
sandpile is stable if its slope does not exceed 1,
i.e., if $\Lip[u(\cdot, t)] \leq 1$. The condition (\ref{UnstLip})
implies that the initial profile $g$ is unstable.
According to the model above, the initially unstable sandpile collapses
instantaneously to
a stable one, with the height function $u$.

The transformation $g \rightarrow u$ was studied in  \cite{EFG}
by relating it to the
evolution problem governed by Monge-Kantorovich mass transfer.
Introduce the convex functional 
$I_{\infty}: L^2(\bR^n) \rightarrow \bR^1 \cup \{ \infty \}$ defined by
$$
 I_{\infty} [ v ] =
\left\{  
\begin{array}{ll}
0 & \mbox{if $v \in L^2(\bR^n),\;\;\; \mid Dv \mid \leq 1\; $ a.e.} \\
+\infty& \mbox{otherwise}.
\end{array}
\right.
$$
Then 
$u(x)=w(x,1)$, where $w$ is the unique solution of the evolution problem
\begin{equation}
\label{evW}
\left\{  
\begin{array}{ll}
\frac{w}{t} - \partial_t  w
 \in  \partial I_{\infty} [ w ],  & \mbox{for a.e. } t \in [\frac{1}{L}, 1], \\
w=\frac{g}{L}=\hat{g}, & \mbox{at } t=\frac{1}{L},
\end{array}
\right.
\end{equation}
with $L=\mbox{Lip} [g] > 1$. 
Note that $\mbox{Lip}[\hat{g}] =1$, and thus the evolution problem
(\ref{evW})
is well-posed. The connection between this evolution problem and
Monge-Kantorovich mass transfer is explained in \cite{EFG}.

In the static case the limits as $p \rightarrow \infty$ of 
the boundary value problems for $p$-Laplacian type equations
were considered in the papers \cite{BDM}, \cite{Jan}
and references therein. 
In many cases the limits
have the form of the distance function to the boundary of the domain.

A limit evolution problem that represents a stable case of
the problem (\ref{pHeatUnS}) was considered in the paper 
\cite{AEW} as a model of a growing stable
sandpile. Let $u(x,t)$ be, as before, the height of the pile,
and let $f(x,t)\geq 0$ be the sand source function which describes
the rate at which sand is coming at location $x$ at time $t$.
The dynamics is
defined as following. Consider the nonhomogeneous evolution problem
\begin{equation}\label{pHeatStabl}
\left\{  
\begin{array}{ll}
\partial_t u_{p} - \mbox{div}(\mid Du_p \mid ^ {p-2} Du_p ) =f, & 
\mbox{in}\;\; \bR^n \times (0, \infty), \\
u_p=0, & \mbox{on}\;\; \bR^n \times \{0\},
\end{array}
\right.
\end{equation}
and let $p \rightarrow \infty$. Then for any $T>0$
$$
u_p \rightarrow u \;\;\;\;\;\;\;\; \mbox{uniformly on} \; \bR^n \times [0, T]
$$
and this limit $u$ satisfies the equation
\begin{equation}
\label{evU}
f - \partial_t  u
 \in  \partial I_{\infty} [ u ],  \mbox{for a.e. } t \geq 0
\end{equation}
with zero initial condition.
In particular the case was considered when sand sources are concentrated 
in $m$ points,
i.e., if
$$
f(x,t)= \sum_{k=1}^m f_k(t) \delta_{d_k}(x),
$$
where $ \delta_{d_k}(\cdot)$ 
is the Dirac mass at the point $d_k \in \bR^n$. The
solution has the explicit form
$$
u=\max(0, z_1(t)-\mid x-d_1 \mid, ..., z_m(t)-\mid x-d_m \mid),
$$
where the height functions $\{ z_k(t) \}_{k=1}^m$ satisfy a system 
of ODE. Thus the
solution has the form of growing and interacting cones centered at the 
points $\{d_k \}_{k=1}^m $.

By analogy with the static problems and the model of growing sandpiles
\cite{AEW}, one can expect to find solutions of (\ref{evW}) that
have the form of distance function to the moving boundary or 
superposition of such, i.e., growing and interacting cones. 
However, unlike
the situation considered in \cite{AEW}, interacting spherical cones
do not maintain their shape under the evolution defined by  (\ref{evW}).
A more suitable class of solutions is the following class of more 
general ``cones''.  A single cone solution has the form
\begin{equation}
\label{2}
w(x,t) = d_{\Omega_t}(x),
\end{equation}
where $\{ \Omega_t \}$ is a suitable expanding family of sets.
A multiple interacting cones solution
has the form (for simplicity we consider the case of two cones)
\begin{equation}
\label{2'}
w(x,t)=\max(d_{\Omega^1_t}(x),d_{\Omega^2_t}(x)),
\end{equation}
where $\{ \Omega_t^1 \}$, $\{ \Omega_t^2 \}$
are suitable expanding families of sets.
Thus in order to find such solutions one has to find the appropriate
families of sets.

The following equations of motion of  the 
boundaries of the sets in (\ref{2}) and (\ref{2'})  are 
derived heuristically
in \cite{EFG}. 
Assume that ${\bf \Gamma}$ defined by (\ref{defE})
 is a smooth surface in $\bR^n \times \bR^1$. 
Let $y \in \Gamma_t$ and let $\kappa_1, ..., \kappa_{n-1}$
be the principal curvatures of $\Gamma_t$ at $y$.
Denote $\kappa = (\kappa_1, ..., \kappa_{n-1})$.
Let $\gamma(y)=\gamma_{\Omega_t}(y)$ be the radius of the largest ball 
touching $\Gamma_t$ at $y$
from within $\Omega_t$. Let $V(y)$ be the outer normal
velocity of  $\Gamma _t$ at the point $y$
defined by (\ref{outNormVelDef}).
Then in the case of a single cone solution (\ref{2})
the equation for the surface $\Gamma_t$ is:
\begin{eqnarray}
\label{1}
 V(y)\;\; & = & \frac{1}{t}
F(\kappa(y), \gamma(y))\;\; \mbox{ for } y \in \Gamma_t, 
 \nonumber \\
\mbox{where}\;\;\;\;\;\; & & \\
 F(\kappa, \gamma) & = &
\frac{\int^{\gamma }_0 s\prod^{n-1}_{i=1} (1-\kappa_i
s)ds}{\int^{\gamma }_0 \prod^{n-1}_{i=1} (1-\kappa_i
s)ds}. \nonumber
\end{eqnarray}

To write equations for the motion of the surfaces $\Gamma_t^k =
\partial\Omega_t^k,\; k=1,2$ in the case of the interacting cones
solution (\ref{2'}), we extend the function
$\gamma(\cdot)=\gamma_{\Omega}(\cdot)$ of a set $\Omega \in \bR^n$
to a function $\bR^n\rightarrow \bR^1$ by setting
$\gamma \equiv 0$ on $\bR^n\backslash \overline{\Omega}$, and
$\gamma (x)=\gamma(y)$ for $x \in \Omega$ such that
$y \in {\cal N}_{\partial \Omega}(x)$. Then the equations are:
\begin{eqnarray}
\label{1'}
V(y)\;\;\;\;\;\;\;
  & = & \frac{1}{t} F(\kappa(y), \gamma_{\Omega_t^k}(y),  
\gamma_{\Omega_t^1 \cap \Omega_t^2}(y) )\;\; \mbox{ for } y \in \Gamma_t^k,
\;k=1,2,   
 \nonumber \\
\mbox{where}\;\;\;\;\;\;\;\;\;\;\; & & \\
F(\kappa, \gamma, \delta) & = &
\frac{\int^{\gamma}_{\delta} s\prod^{n-1}_{i=1} (1-{\kappa}_i
s)ds}{\int^{\gamma}_{\delta} \prod^{n-1}_{i=1}
(1-{\kappa}_i
s)ds}.
 \nonumber
\end{eqnarray}

 The equations (\ref{1}) and (\ref{1'}) have
the form of nonlocal geometric curvature motion.
Geometric nonlocality
is caused by the functions $\gamma(\cdot)$.

In this paper we prove that if the moving
surfaces satisfy geometric equations (\ref{1}) or (\ref{1'}) and
some regularity conditions, then the functions
(\ref{2}) or (\ref{2'}) respectively are solutions of the 
evolution problem (\ref{evW}). Namely, we
 prove the following.

\begin{theorem}
\label{Theorem 1}
  Let $\{\Omega_t \}$  be a locally
Lipschitz continuous
family of open convex bounded sets,
where $ t\in [a,b]$, $ b>a>0$. Suppose that the equation (\ref{1}) 
is satisfied at every point  of $(2,1)$ differentiability of the surface
${\bf \Gamma}$ defined by (\ref{defE}). 
Then the function  
(\ref{2})
is a solution of the evolution equation
\begin{equation}
\label{3}
{w\over t}-\partial_t w\in \partial I_\infty [w]\;\;\; \; \mbox{for
}a.e.\;\;t\in [a,b].
\end{equation}
\end{theorem}

\begin{theorem}
\label{Theorem 2}
Let $\{ \Omega^1_t \}$, $\{ \Omega^2_t \}$
 be
two Lipschitz continuous families of 
open bounded convex sets, where $ t\in
[a,b]$,  $b>a>0$.
Suppose that the equation (\ref{1'}) is satisfied
at every
point  of $(2,1)$ differentiability of the surfaces
${\bf \Gamma}^k = \cup_t(\Gamma^k_t \times \{ t \}), \;k=1,2$.
Then the function (\ref{2'})
is a solution of
 the evolution equation (\ref{3}).
\end{theorem}

The convexity condition in Theorems \ref{Theorem 1}, \ref{Theorem 2} can
be relaxed to the following "semiconvexity" condition.
\begin{condition}
\label{Cond.1}
{ \bf (lower curvature bound condition with radius $r$ ).} 
An open set $\Omega \subset \bR^n$ satisfies this condition 
if  $\gamma_{ \bR^n  \setminus \overline{\Omega}}(x) \geq 2r$
for any $x \in  \bR^n  \setminus \overline{\Omega }$.
\end{condition}
It is easy to check  that if a set is either
convex or has $C^{1,1}$ boundary,
then the set satisfies the Condition \ref{Cond.1}. 
\begin{remark}
\label{remrkCond 1}
Let the set $\Omega$ satisfy Condition \ref{Cond.1}. Let 
$\partial \Omega$ be twice differentiable at a point $y$,
i.e., (\ref{SurfDiff}) and (\ref{SurfTwiceDiff})
are satisfied in an appropriate coordinate system on $\bR^n$. 
Then the principal curvatures of $\partial \Omega$ at
$y$ satisfy
$$
\kappa_i \geq -\frac{1}{2r} \;\;\; \;
\mbox{for } \; i=1,...,n-1.
$$
 This follows from
the fact that
the graph of the function $\Psi$ defined by 
(\ref{SurfDiff}) and (\ref{SurfTwiceDiff}) lies above the
ball $B_{2r}^n(0, -2r)$ by  Condition
\ref{Cond.1}.
\end{remark}

The following two theorems state that convexity can be replaced by the
weaker Condition \ref{Cond.1} in Theorems \ref{Theorem 1} and
\ref{Theorem 2}.
\begin{theorem}
\label{Theorem 3}
  Let $\{\Omega_t \}$  be a locally
Lipschitz continuous
family of open bounded sets,
where $ t\in [a,b]$, $ b>a>0$.
Let
for every $ t\in [a,b]$ the set $\Omega_t$ satisfy Condition
\ref{Cond.1} with radius $r_0 > 0$.
 Suppose that the equation (\ref{1}) 
is satisfied at every point  of $(2,1)$ differentiability of the surface
${\bf \Gamma}$ defined by (\ref{defE}). Then the function  
(\ref{2})
is a solution of the evolution equation (\ref{3}).
\end{theorem}

\begin{theorem}
\label{Theorem 4}
Let $\{ \Omega^1_t \}$, $\{ \Omega^2_t \}$
 be
two Lipschitz continuous families of 
open bounded sets, where $ t\in
[a,b]$,  $b>a>0$. Let
for every $ t\in [a,b]$ the sets $\Omega_t^1$, $\Omega_t^2$ satisfy 
Condition \ref{Cond.1}
with radius $r_0 > 0$.
Suppose that the equation (\ref{1'}) is satisfied
at every
point  of $(2,1)$ differentiability of the surfaces
${\bf \Gamma}^k = \cup_t(\Gamma^k_t \times \{ t \}), \;k=1,2$.
Then the function (\ref{2'})
is a solution of
 the evolution equation (\ref{3}).
\end{theorem}

Theorems \ref{Theorem 1} and \ref{Theorem 2} follow from
Theorems \ref{Theorem 3} and \ref{Theorem 4} since convex sets
satisfy Condition \ref{Cond.1} with any radius.

The proof of Theorems \ref{Theorem 3} and \ref{Theorem 4} utilizes the
following relation of the equation
(\ref{3})  and the Monge-Kantorovich mass transfer problem. 
Fix $t$.
Equation 
(\ref{3})  
implies that that the function $w(\cdot,t)$ is the Monge's 
potential for the optimal
transfer
of the measure with density $\frac{w}{t}(\cdot,t)$  
into the measure with density
$\partial_t{w}(\cdot,t)$.  
Then, according to \cite{EG}, there exists a measurable
 function $a(x)$, the {\it mass
transport density}, 
that satisfies the following properties:
\begin{eqnarray}
a \geq 0, \;\;\; & \; & \mbox{supp}(a) \subset \{x \;\mid\;\mid Dw\mid=1\}, 
\nonumber \\
 & & \label{defMassTrDen}  \\
-\mbox{div}(aDw) & = & \frac{w}{t} - \partial_tw, \nonumber
\end{eqnarray}
where the last equation is understood in the weak sense and $t$ is fixed. 
Conversely,
if there exists a function $a(x)$ satisfying (\ref{defMassTrDen}),
then the function $w$ is a solution of (\ref{3}) for
the given $t$, see Section \ref{Th.3} below.

The idea of the proof of Theorem \ref{Theorem 3} is as follows.
Let a family of sets $\{\Omega_t\}$ satisfy (\ref{2}).
An open bounded set can be represented as the union of distance rays,
i.e., maximal intervals that start at the boundary, on which distance to the 
boundary is the linear function with slope 1. The Condition \ref{Cond.1}
implies that the collection of distance rays possesses "nice"
measure-theoretic properties. Fix $t$ and consider a distance ray $R_y$ of 
$\Omega_t$ starting 
at a smooth point $y\in \partial \Omega_t$ . Introduce the coordinate 
$s$ on $R_y$, the 
distance to the boundary.
 The equation (\ref{defMassTrDen}) 
can be formally rewritten as:
\begin{equation}
\label{frmlEq}
Da \cdot Dw + a \Delta w - \frac{w}{t} -  \partial_tw = 0.
\end{equation}
Since $w(\cdot, t) = \dist(\cdot, \Gamma_t)$, we have on $R_y$
$$
Da \cdot Dw = \frac{da}{ds}, \;\;\; w=s ,\;\;\;
-\Delta w = \sum^{n-1}_{i=1} 
{{\kappa_i}\over
{1-\kappa_i s}},  \;\;\; 
\partial_tw = V(y),
$$
where $\kappa_1, ..., \kappa_{n-1}$ are the principal curvatures 
of $\Gamma_t$ at $y$, and
$V(y)$ is the outer normal velocity of $\Gamma_t$ at $y$.
Thus the equation (\ref{frmlEq})
can be formally rewritten
on the
ray $R_y$ as the following ODE
\begin{equation}
\label{TrDenODE}
a'(s)-a(s)\sum^{n-1}_{i=1} 
{{\kappa_i}\over
{1-\kappa_i s}} -{s \over t}+V(y)=0,
\end{equation}
This ODE
has a solution $a(s)$ with zero boundary conditions at both ends of the
ray $R_y$ if $V(y)$ satisfies (\ref{1}). Define $a(x,t)$ by ODE
(\ref{TrDenODE}) with zero boundary conditions
on each distance ray that starts
at a smooth point of the boundary. We prove that $a(\cdot, t)$ is 
a measurable function that satisfies (\ref{defMassTrDen}). In order to
do this we examine the properties of distance function and distance rays,
and use a nonsmooth (Lipschitz) change of coordinates on the set 
$\Omega_t$ with $n-1$ coordinates along the boundary
and the $n$-th coordinate along distance rays. 

We prove Theorems \ref{Theorem 3} and  \ref{Theorem 4} 
in Sections \ref{Sets} - \ref{Th.4}.
In Section \ref{Sets} we examine properties of sets 
satisfying Condition \ref{Cond.1}, and in particular
we describe the change of coordinates mentioned above.
In Section  \ref{SectDens}  we define the mass transport density function
and examine some properties of this function. In Section \ref{ChngVar}  
we prove that the mass balance equation (\ref{1}) is satisfied. 
In Section \ref{Th.3} we conclude the
proof of Theorem \ref{Theorem 3}. Finally, in Section \ref{Th.4} we 
sketch  the proof of Theorem \ref{Theorem 4}.

\section{Properties of distance function and ridge sets of the sets satisfying
the lower curvature bound condition}
\label{Sets}

Let $\Omega$ be an open set.
When there is no possible confusion, we write $d(x)$
for $d_{\Omega}(x)$ defined by (\ref{distFunc00}). 

Let $x\in \overline{\Omega}$. 
Denote by $R_x$ the longest line segment through $x$ in
$\overline{\Omega}$ along
which $d_{\Omega}(\cdot)$ is a linear function with slope 1. 
We call $R_x$ the distance
ray of $x$. Note that a point $x$ can have more
than one distance ray, in such case denote any of them by $R_x$. 
If $x \in \Omega$ then all such rays $R_x$ have equal length.
If $x \in \partial \Omega$ then rays $R_x$ in the case of
general open set $\Omega$ can
have different length. However if $\Omega$ satisfies
 Condition \ref{Cond.1}
then there exists at most one $R_x$ for $x \in 
\partial \Omega$, see Proposition \ref{Proposition 2} below.

For any $x \in \Omega$ one endpoint
of $R_x$ lies on $\partial \Omega$ and belongs to
the set ${\cal N}_{\partial \Omega}(x)$. Call this endpoint the
lower end of $R_x$.  Call another endpoint of
$R_x$ the upper end of $R_x$ (the names ``lower'' and ``upper''
correspond to the positions of the endpoints of $R_x$ on
the graph of $\;y=d_{\Omega}(x)$). Let $y, \;\; z$ be upper
and lower ends of the ray $R$, and let $x=\lambda y + 
(1-\lambda) z$ where $\lambda \in (0, \;1)$. Then we say that
the point $x$ lies in the relative interior of the ray $R$.

\begin{remark}
\label{remUnNr}
If $x \in \Omega$ lies in the relative interior of a
distance ray $R_x$,
then the ray $R_x$ is the unique distance ray that contains
$x$, and the set  ${\cal N}_{\partial \Omega}(x)$
consists of one point, see \cite{EH}.
\end{remark}

\begin{remark}
\label{remPosReach}
It follows from Remark \ref{remUnNr} that sets satisfying
Condition \ref{Cond.1}
are sets with positive reach defined by H. Federer \cite{FedCMeas}.
More precisely, if a set
$\Omega$ satisfies  Condition \ref{Cond.1} with
radius $r$ then reach($\Omega$)=$2r$.
\end{remark}

The function $\gamma_{\Omega}(\cdot)$ in the equations
(\ref{1}), (\ref{1'}) is called ridge
function and can be expressed as following:
\begin{equation}
\label{defGamma}
 \gamma_{\Omega}(x) =
\left\{  
\begin{array}{ll}
0 & x \in \bR^n \setminus \overline{\Omega}, \\
\mid R_x \mid & x \in \Omega, \\
\sup \mid R_x \mid  & x \in \partial \Omega.
\end{array}
\right.
\end{equation}
The function $\gamma_{\Omega}(\cdot)$ is uppersemicontinuous
in $\bR^n$, cf. \cite{EH}, Proposition 3.2.

Define the set
$$
{\cal R} = \{\; x \in \overline{\Omega} \;\;\; | \;\;\;
d_{\Omega}(x) =  \gamma_{\Omega}(x)\; \}.
$$
The set ${\cal R}$ is called {\it ridge set} of $\Omega$. 
Note that if  $x \in {\cal R}$ then either $x$ is the upper end 
of some
distance ray of $\Omega$,
or otherwise $x$ lies at  $ \partial \Omega$ and there are no
points $y\in \Omega$ such that 
$x \in {\cal N}_{\partial \Omega}(y)$.

\begin{proposition}
\label{Proposition 2} 
Let $\Omega \subset \bR^n$ be a bounded open set satisfying 
Condition \ref{Cond.1} with radius $r>0$. 
Let $\Omega_r$ be the $r$-neighborhood
of $\Omega.$  Then 
  
(a) $\gamma_{\bR^n \setminus \overline{\Omega}_r }(x) \geq r$ 
for all $x \in \bR^n \setminus \Omega_r$, i.e., the 
set $\Omega_r$ satisfies 
Condition \ref{Cond.1} with radius $r\over 2$.
  
(b) For any $y_r \in \partial \Omega_r$ there exists a unique $y \in 
\partial \Omega$ such that $|y_r-y|=r$, and $|y_r-x|>r$
for all $x\in \overline{\Omega}\backslash \{ y \}$.  For any 
$y\in \partial \Omega$ 
there exists 
$y_r \in \partial \Omega_r$ 
such that
$|y-y_r|=r$.  If $\partial \Omega$ is differentiable at $y$
  then
 $y-y_r$ is orthogonal to $\partial \Omega$ at $y$.  If there exists a ball
inside $\Omega$ touching $\partial \Omega$ at $y$, then $\partial \Omega$
is differentiable at $y$. In particular, $\partial \Omega_r$ is everywhere 
differentiable.
    
(c) For any $x\in \Omega$ we have
\begin{equation}
\label{distRelat}
\mbox{dist}(x,\partial \Omega_r)=\mbox{dist}(x,\partial \Omega)+r.
\end{equation}
If  $x \in {\Omega_r} \setminus  \overline{\Omega}$ then 
 $\dist(x, \partial \Omega_r) < r$ and $x \notin {\cal R}_r$, 
where ${\cal R}_r $ is ridge set of $\Omega_r$. 
    
(d) Sets $\Omega$ and $\Omega_r$ have the same ridge set, i.e., 
${\cal R}={\cal R}_r$.
\end{proposition}

\begin{proof} 

(a) Let $x \in \bR^n \setminus \overline{\Omega}_r$. 
Then $\dist(x, \Omega) = R > r$. Let 
$y \in  \partial{\Omega}$
be a point such that $\mid x-y \mid = R$. Let $y_r$ be such point of 
the interval
connecting $x$ and $y$ that $\mid y-y_r \mid = r$. Then 
$y_r \in \overline{\Omega}_r$.

Consider first the case $R \geq 2r$.
Then $\dist(x, \Omega_r)\geq r$. Since if this is not true then there exists
such point  $z_1 \in \overline{\Omega}_r$ that $\mid x-z_1 \mid < r$,
and there exists $z_2 \in \overline\Omega $ such that 
$\mid z_1 - z_2 \mid \leq r$. Now we get 
$R \leq \mid x-z_2\mid \leq \mid x-z_1 \mid + \mid z_1 - z_2 \mid < 2r$,
a contradiction. So  $\gamma_{\bR^n \setminus \overline{\Omega}_r }(x)
 \geq \dist(x, \Omega_r)\geq r$.

 Consider now the case $r<R < 2r$. Then by Condition \ref{Cond.1} we have
$ \gamma_{\bR^n \setminus 
\overline{\Omega} }(x) \geq 2r$ and so
$R=\dist(x, \Omega) < \gamma_{\bR^n \setminus 
\overline{\Omega} }(x)$. Thus $y$ is the unique 
point on $\partial \Omega$ such that $\mid x-y \mid = \dist(x,\Omega)$.
Let $z = y + 2r\frac{x-y}{\mid x-y \mid}$. Then from the condition  
\ref{Cond.1}
and \cite{EH}, Lemma 3.4, it follows that for any $\tilde{x}$
of the form $tx+(1-t)z$ 
where $0<t<1$
the point $y$ is the unique point on $\partial \Omega$ such that  
$\dist(\tilde{x}, \partial \Omega)
= \mid \tilde{x} - y \mid$. Then for each such point $\tilde{x}$, the 
point $y_r$ is 
 the unique point on $\partial \Omega_r$ such that  $\dist(\tilde{x}, 
\partial \Omega_r)
= \mid \tilde{x} - y_r \mid$. For suppose on the
 contrary that there exists $z_1 \in 
\partial \Omega_r$
such that $z_1 \neq y_r$ and $\mid \tilde{x} - z_1 \mid \leq \mid 
\tilde{x} - y_r\mid $.
Then there exists $z_2 \in \partial \Omega$ such that 
$\mid z_1 - z_2  \mid \leq r$. Consider first the case $z_2 \neq y$. Then
$$
\mid \tilde{x} - y \mid = \mid \tilde{x} - y_r \mid + \mid y_r - y \mid =
\mid \tilde{x} - y_r \mid + r  \geq \mid \tilde{x} - z_1 \mid +
\mid z_1 - z_2 \mid \geq \mid \tilde{x} - z_2 \mid,
$$
a contradiction with the fact that $y$ is the unique nearest to $\tilde{x}$
point
on $\partial \Omega$. Consider now the case $z_2 = y$. Since 
the ball of radius $R$ and center  $\tilde{x}$ lies outside $\Omega$
and $|z_1-y| \leq r$, it
follows that $z_1$ does not lie on the interval connecting
$\tilde{x}$ and $ y$, and so  
$$
\mid \tilde{x} - y \mid < \mid \tilde{x} - z_1 \mid + \mid z_1 - y \mid \leq
\mid \tilde{x} - y_r \mid +r = \mid \tilde{x} - y \mid,
$$
 a contradiction.
So, $y_r$ is the unique nearest to  $\tilde{x}$ point on $\partial \Omega_r$
for any $\tilde{x}$ on the interval
connecting $z$ and $y_r$, so
$\gamma_{\bR^n \setminus \overline{\Omega}_r}(x) \geq \mid z - y_r \mid  = r$.
  \\

(b) Let $y_r \in \partial \Omega_r$. Then $\dist(y_r, \Omega) = r$. 
Thus there exists $y\in \partial \Omega$ such that  $\mid y-y_r \mid = r$.
By Condition \ref{Cond.1} , $\gamma_{\bR^n \setminus 
\overline{\Omega}}(y_r) \geq 2r
> \dist(y_r, \partial \Omega)$, and thus by \cite{EH}, Lemma 3.4,
$y$ is the unique point of $\partial \Omega$ nearest to $y_r$.

Let $y\in \partial \Omega$. Let $y_i \in \bR^n \setminus 
\overline{\Omega} $ for $i=1,2,...$, and 
let $y_i \rightarrow y$. By condition \ref{Cond.1}, 
$\gamma_{\bR^n \setminus 
\overline{\Omega}}(y_i) \geq 2r$, and by uppersemicontinuity of
$\gamma_{\bR^n \setminus \overline{\Omega}}(\cdot)$ we get 
 $\gamma_{\bR^n \setminus \overline{\Omega}}(y) \geq 2r$.
Thus we have shown
that every point of $\partial  \Omega$  has an exterior tangent ball
of the radius $2r$. So if  a point $y$ of $\partial  \Omega$ has 
an interior tangent ball, then $\partial  \Omega$ is differentiable
at $y$.

Let $y \in \partial  \Omega$. Then, as we have shown, 
$\gamma_{\bR^n \setminus \overline{\Omega}}(y) \geq 2r$.
Thus there exists a distance ray $R_y$ of the set $\bR^n 
\setminus \overline{\Omega}$ starting at $y$ such
that 
the length of $R_y$ is at least
$2r$.  Let $e$ 
be the unit vector in the direction of $R_y$, let
$y_r = y + re$. Then ${\cal N}_{\partial \Omega}(y_r) = \{ y \}$
by Remark \ref{remUnNr}.
Since $\mid y-y_r \mid = r$ we get $y_r \in \partial \Omega_r$.

We have shown that the interval connecting $y$ and $y_r$ lies on the
distance ray of the set $\bR^n \setminus \overline{\Omega}$ starting at $y$,
and on distance ray of the set $\Omega_r$ starting at $y_r$.  It follows that
the vector $y-y_r$ is orthogonal to $\partial \Omega_r$ at $y_r$  and also
to  $\partial \Omega$ at $y$ if $\partial \Omega$ is differentiable at $y$.
  \\

(c) Let $x\in \Omega^0, \, y\in \partial \Omega$, and $|x-y|=
\dist(x, \partial \Omega)$. 
Then, by (b),  
$$
y_r := y + \frac{r}{|y-x|}(y-x) \in \partial \Omega_r.
$$
Clearly,
$$
|x-y_r| = |x-y| + r = d_{\Omega}(x) + r.
$$
Assume that there exists $z_r \in \partial \Omega_r$ such
that 
$$
|x - z_r| < d_{\Omega}(x) + r.
$$ 
Let $z$ be the point
of intersection of $\partial \Omega$ with the interval 
connecting $z_r$ and $x$ (such point exists since
$x \in \Omega, z_r \notin \overline{\Omega}$).
Then $|z-z_r| \geq r$ and so
$$
|x - z| < d_{\Omega}(x).
$$
This contradicts the fact that $z \in \partial \Omega$.
Thus $y_r \in {\cal N}_{\partial \Omega_r}(x)$
and thus (\ref{distRelat}) is proved.

It remains to prove the last assertion of the statement (c).

Let $x\in {\Omega}_r\backslash 
\overline{\Omega}$ . Let  $y_r \in {\cal N}_{\partial \Omega_r}(x)$. 
Let $R_{y_r}$ be the ray orthogonal to
$\partial \Omega_r$ at $y_r$ and let $y=R_{y_r} \cap \partial \Omega$. 
Then   by (b) we have the following:
$x\in R_{y_r},\;\;x$ lies on the interval connecting $y$ and $y_r$, and $|y-
y_r|=r=\mbox{dist}(y,\partial \Omega_r)$.  Then  
it follows from \cite{EH}, Lemma 3.4 that  $x \notin {\cal R}_r$. Also,
$\dist(x, \partial \Omega_r) = \mid x-y_r \mid < |y-
y_r|=r$.
  \\

(d)  We show first that ${\cal R}\subset {\cal R}_r$.
 If $x\in \Omega$ and $B_\rho(x)\subset \Omega$, then
$B_{\rho +r}(x)\subset \overline{\Omega}_r$.
Let $x\in {\cal R}\backslash \partial \Omega$. Let 
$y \in {\cal N}_{\partial \Omega}(x)$ and let
$y_r=y+r\frac{y-x}{|y-x|}$.  Then by (b),(c) we see that 
$y_r \in \partial \Omega_r,$ and $|x-y_r|=\dist(x, \partial \Omega) +r=
\mbox{dist}(x,\partial
\Omega_r)$. 
Let $\tilde{x}=y+(1+\varepsilon)(x-y)$ where $\varepsilon>0$ is 
small enough so that 
$\tilde{x}\in \Omega$.  Since $x\in {\cal R}$, it follows that 
 $\dist(\tilde{x}, \partial \Omega)<|\tilde{x}-
y|=(1+\varepsilon)\dist(x, \partial \Omega)$.  Let 
$\tilde{y}\in\partial\Omega$ be such
point that $\dist(\tilde{x}, \partial \Omega)=|\tilde{x}-\tilde{y}|$, let
$\tilde{y}_r=\tilde{y} +r\frac{\tilde{y}-
\tilde{x}}{|\tilde{y}-\tilde{x}|}$, then $\tilde{y}_r\in
\partial \Omega_r$ and $|\tilde{x}-
\tilde{y}_r|=\dist(\tilde{x}, \partial \Omega_r)$.  Then
$$
\mbox{dist}(\tilde{x}, \partial \Omega_r)=
\mid  \tilde{y}-\tilde{x} \mid + r
<(1+\varepsilon)
\dist(x, \partial \Omega)+r
=|\tilde{x}-y_r|,
$$
so $\gamma_{\Omega_r}(x)= |x-y_r| = d_{\Omega_r}(x)$. Thus $x\in {\cal R}_r$.

Let $x\in{\cal R}\cap \partial \Omega$,  i.e., there does not exist a ball
inside $\Omega$ that touches $\partial \Omega$ at $x$.  
Let $B_r(y_r)$ be a ball of
radius $r$ and center at such point $y_r$ that 
$\overline{B}_r(y_r)\cap \overline{ \Omega}=\{x\}$,  such ball exists by 
Condition \ref{Cond.1}.  Then  $\dist(y_r, \partial
\Omega)=r$, and by (b) $x$ is the unique point of $\partial \Omega$
such that  $|x - y_r| = r$.
Let $\tilde{x}=y_r+(1+\varepsilon)(x-y_r)$ for 
$\varepsilon\in(0, {r \over 2})$. 
Given $\varepsilon$, there are two possibilities:
either  $\tilde{x} \in \Omega$ or $\tilde{x} \notin \Omega$.

Let $\tilde{x} \in \Omega$.  Then, since $x\in {\cal R}\cap \partial
\Omega,\;\;\dist(\tilde{x}, \partial \Omega)<|x-\tilde{x}|=\varepsilon r$.
Let $\tilde{y}\in \partial \Omega$ be such that $|\tilde{x}-
\tilde{y}|=\dist(\tilde{x}, \partial \Omega)$, let $\tilde{y}_r=\tilde{y}+
r\frac{\tilde{y}-
\tilde{x}}{|\tilde{y}-\tilde{x}|}$. Then $\tilde{y}_r \in \partial
\Omega_r$ and by (c) we get
$|\tilde{x}-\tilde{y}_r|=\dist(\tilde{x}, \partial \Omega)+r=
\mbox{dist}
(\tilde{x},\partial \Omega_r)$, and so $\mbox{dist}(\tilde{x}, \partial
\Omega_r)<(\varepsilon +1)r= |\tilde{x}-y_r| $, so 
$y_r \notin {\cal N}_{\partial \Omega_r}(\tilde{x})$.

Let now  $\tilde{x} \notin \Omega$, i.e., $\tilde{x} \in     
\Omega_r \setminus \overline{\Omega}$ since $\varepsilon < {r \over 2}$.
 Then, by (c),  $\dist(\tilde{x}, \partial \Omega_r) < r$,
but $|\tilde{x}-y_r| = (\varepsilon +1)r$, so 
$y_r \notin {\cal N}_{\partial \Omega_r}(\tilde{x})$. 

Thus
$y_r \notin {\cal N}_{\partial \Omega_r}(\tilde{x})$ for
any small $\varepsilon>0$. This implies $x \in {\cal R}_r.$

Thus ${\cal R}\subset {\cal R}_r.$ 

Now we show that ${\cal R}_r \subset {\cal R} $. 
By (b), (c) we get 
${\cal R}_r \cap (\overline{\Omega}_r \setminus  \overline{\Omega} ) = 
\emptyset$.

It remains to consider $x\in \overline{\Omega}\cap {\cal R}_r $.  Let 
first $x\in \Omega
\cap {\cal R}_r$.  Let as above $y \in {\cal N}_{\partial \Omega}(x)$,
$ y_r =y+r\frac{y-x}{|y-x|} \in {\cal N}_{\partial \Omega_r}(x)$,
$\tilde{x}=y+(1+\varepsilon)(x-y)$ where $\varepsilon>0$ is 
small. 
From $x \in {\cal R}_r$ we conclude that 
 $\mbox{dist}(\tilde{x}, \partial \Omega_r)<|\tilde{x}-
y_r|=(1+\varepsilon)\dist(x, \partial \Omega)+r$.  Let $\tilde{y}_r\in 
\partial
\Omega_r$ be such that $|\tilde{x}-\tilde{y}_r|=\mbox{dist}(\tilde{x}, 
\partial
\Omega_r)$ and let $\tilde{y}$ be the unique point of $\partial \Omega$ 
nearest to
$\tilde{y}_r$, then $|\tilde{y}-
\tilde{y}_r|=r$ (uniqueness of $\tilde{y}$ follows
from (b)).  Then by (c), (b)
the points $\tilde{x},\tilde{y},\tilde{y}_r$ lie on one distance ray
 and $\dist(\tilde{x}, \partial \Omega)=
|\tilde{x}-\tilde{y}|$, so
$$
\dist(\tilde{x}, \partial \Omega)=|\tilde{x}-\tilde{y}|=|\tilde{y}_r-
\tilde{x}|-
|\tilde{y}_r-y|<|\tilde{x}-y_r|-r=
$$
$$
(1+\varepsilon)\dist(x, \partial \Omega)=(1+\varepsilon)|x-
y|=|\tilde{x}-y|
$$ 
and so $\tilde{x}\in {\cal R}$. \\
The case of $x\in \partial \Omega \cap {\cal R}_r$ is similar to the
case above. 
\end{proof}

\begin{remark}
\label{RemSet1}
Examples of nonconvex polygons on a plane show that
 without assuming some condition of the type
of Condition \ref{Cond.1} the 
statements (b)-(d) of Proposition \ref{Proposition 2} are not true in general.
\end{remark}

\begin{proposition}
\label{Proposition 3}
Let $\Omega \in R^n$ satisfy Condition \ref{Cond.1} with radius $r$. Then 
$\partial \Omega_{\rho}\in C^{1,1}$ for any $\rho \in (0, \frac{r}{M})$.
\end{proposition}
\begin{proof}
The assertion follows from inequality 4.8(8) of \cite{FedCMeas}
 and Proposition \ref{Proposition 2} (d). Also, inequality (\ref{1.1})
proved in Appendix A1 can be applied instead of  4.8(8) of \cite{FedCMeas}.
\end{proof}

Recall the following facts (see e.g. \cite{GilTr}, Lemma 14.16, 14.17).
Let $\Omega \subset \bR^n $ be an open bounded set with $C^2$ boundary
$\partial \Omega$. 
Let $\Omega_r$ be, as above, $r$-neighborhood of $\Omega$.
Then $\partial \Omega_r\in C^2$ for small $r$. If 
$y\in \partial \Omega$ and if $y_r\in \partial
\Omega_r$ is the unique point on $\partial \Omega_r$ such that $|y-
y_r|=r$ then the principal coordinate systems of $\partial \Omega$ at
$y$ and of $\partial \Omega_r$ at $y_r$ are parallel. Denote by
$\kappa_i$ and  $\kappa_{i,r},\;\; (i=1,...n-1)$ the principal 
curvatures of $\partial \Omega$ at $y$ and of $\partial \Omega_r$ 
at $y_r$ respectively. We have
\begin{equation}
\label{relNbhd}
\kappa_{i,r}=\frac{\kappa_i}{1+\kappa_ir}\;\;\;
\mbox{ for }\;\;i=1,...,n-1.
\end{equation}
In the next two propositions we
show similar facts for sets with nonsmooth boundaries satisfying Condition
\ref{Cond.1}
and for Lipschitz families of such sets.

\begin{proposition}
\label{Proposition 4}
Let $\Omega \subset \bR^n$ be an open set satisfying Condition
\ref{Cond.1} with radius $r_0$.  Let $0<r \leq r_0$. Fix
$y\in \partial \Omega$. Let $y_r\in \partial
\Omega_r$ be such point that $|y-y_r|=r$.  Then:
  
a) Let $\partial \Omega$ be twice 
differentiable at $y$, i.e.,
in a suitable coordinate system in $\bR^n$ we have
$ y=0,\; \; \Omega
\cap
B_\delta(0)=\{ x_n>\Psi(x')\} \cap
B_\delta(0)$ for small $\delta>0$ and $\Psi(x')=
\sum^{n-1}_{i=1} \kappa_ix^2_i+o(|x'|^2)$.
Then the point $y_r$ is unique given the point $y$,
and in the coordinate system introduced above 
$y_r=(0,-r)$,  $\Omega_r\cap B_{\delta_1}(y_r)=\{ x_n>\Psi_r (x')\} 
\cap B_{\delta_1}(y_r) $ 
for sufficiently small $\delta_1>0$, where
$\Psi_r(x')=-r+\sum^{n-1}_{i=1}\frac{\kappa_i}
{1+\kappa_ir}x^2_i+o(|x'|^2)$, i.e., $\Omega_r$ is
twice differentiable at $y_r$ and (\ref{relNbhd}) holds. 

b) If $\partial \Omega_r$ is twice  
differentiable at $y_r$
and $\gamma_{\Omega}(y) > 0$ then  
$\partial \Omega$ is twice  differentiable at $y$
 and (\ref{relNbhd}) holds (i.e., the  functions 
$\Psi$ and $\Psi_r$ have same expansions
 at 0 as in a)).
\end{proposition}
\begin{proof}

a) Uniqueness of $y_r$ and the fact that $y_r=(0,-r)$
follow from smoothness of
$\partial \Omega$ at $y$. So the interval 
connecting $y$ and $y_r$ lies on $x_n$-axis.
By Propositions \ref{Proposition 2} and \ref{Proposition 3}
the surface $\partial \Omega_r$ is of class $C^{1,1}$,
and the interval 
connecting $y$ and $y_r$ is orthogonal to $\partial \Omega_r$.
So $\Omega_r\cap B_{\delta_1}(y_r)=\{ x_n>\Psi_r (x')\} 
\cap B_{\delta_1}(y_r) $ for some $\delta_1 > 0$, 
where $\Psi_r$ is a $C^{1,1}$ function.

Let $\varepsilon>0$.  There exists $\sigma>0$ such that 
$$
|\Psi(x')-\sum^{n-1}_{i=1}
\kappa_ix^2_i|<\varepsilon|x'|^2 \;\; \mbox{if} \;\; |x'|< \sigma.
$$
Denote $A_\varepsilon=\{(z',z_n) \in \bR^n \;\; |\;\;  z_n\geq \sum^{n-
1}_{i=1}(\kappa_i-\varepsilon)z^2_i\}$.  Then we get:
$$
\partial \Omega \cap B_\sigma (0)\subset
A_\varepsilon\backslash A_{-\varepsilon}.
$$
Since $\gamma_{ \bR^n  \setminus \overline{\Omega}}(y) \geq 2r$,
the point $y_r$ lies in the relative interior of the distance ray $R_y$ of
the set $ \bR^n  \setminus \overline{\Omega}$. Then
 there exists
$\sigma_1>0$ such that for any point $x_r \in
\cap B_{\sigma_1}(y_r)$
we have $ {\cal N}_{\partial \Omega}(x_r) \subset \partial
\Omega \cap B_{\sigma}(0)$. (Proof: suppose this is false, 
then there exist $\sigma_1>0$ and a sequence $x_r^k \rightarrow y_r$
such that for $x^k \in {\cal N}_{\partial \Omega}(x_r^k)$ 
we get $|x^k-y| > \sigma_1$.
Then, passing to a subsequence, we get $x^k \rightarrow x$
where $x\in \partial \Omega$ and $x \neq y$.
By continuity of distance function we get
$x \in {\cal N}_{\partial \Omega}(y_r)$ and
$x \neq y$. By 
Condition \ref{Cond.1} $y_r$ lies in the  relative
interior of the distance ray $R_y$ of the set 
$\bR^n \setminus \overline{\Omega}$. Thus $y_r$ has a unique
nearest point on $\partial \Omega$. A contradiction.).
So 
\begin{equation}
\label{Aepsilon}
x\in A_\varepsilon \backslash A_{-\varepsilon}\;\;\;\mbox{ if }\;\;
 x \in {\cal N}_{\partial \Omega}(x_r)\;\;\; \mbox{ where }\;\;
x_r \in B_{\sigma_1}(y_r).
\end{equation}

Let $N_{\varepsilon}$ be $r$-neighborhood of 
the set $A_{\varepsilon}$.  
By Remark \ref{remrkCond 1} we have $\kappa_i \geq -\frac{1}{2r}$. 
 Let $\varepsilon < \frac{1}{4r}$.
 Then
(by \cite{GilTr} Lemma 14.16, 14.17), $\partial N_\varepsilon$ is
the graph of a $C^\infty$ function $g_\varepsilon$, and
$\frac{\partial^2g_\varepsilon}{\partial x_i\partial
x_j}(0)=\frac{\kappa_i-\varepsilon}{1+(\kappa_i-
\varepsilon)r}\delta_{ij}$.  Let 
$$
\tilde{N}_{\varepsilon,\delta}=\{(z',z_n)\;\; |\;\;
z_n\geq \sum^{n-1}_{i=1}(\frac{\kappa_i-\varepsilon}{1+(\kappa_i-
\varepsilon)r}-\delta)z^2_1-r\}.
$$  
If $\sigma_2>0$ is small enough then 
\begin{equation}
\label{Nvareps}
\tilde{N}_{\pm\varepsilon,-\varepsilon}\cap B_{\sigma_2}(y_r)\subset
N_{\pm\varepsilon}\cap B_{\sigma_2} (y_r) \subset \tilde{N}_{\pm\varepsilon,
\varepsilon} \cap
B_{\sigma_2} (y_r).
\end{equation}
Thus if we choose $\sigma_3$ small, then 
 we get by (\ref{Aepsilon}), (\ref{Nvareps}) 
$$
\partial \Omega_r \cap B_{\sigma_3}(y_r) \subset  \tilde{N}_{\varepsilon,
\varepsilon}\backslash \tilde{N}_{-\varepsilon,-\varepsilon},
$$ 
which means
that for $x_r =(x'_r, \Psi_r(x'_r))\in \partial \Omega_r \cap
B_{\sigma_3}(y_r)$ we have
$$
\sum^{n-1}_{i=1}(\frac{\kappa_i-\varepsilon}{1+(\kappa_i-
\varepsilon)r}-\varepsilon)x^2_{r,i}-r\leq \Psi_r(x'_r)\leq \sum^{n-
1}_{i=1}(\frac{\kappa_i+\varepsilon}{1+(\kappa_i+
\varepsilon)r}+\varepsilon)x^2_{r,i}-r.
$$
Sending $\varepsilon$ to 0 we conclude the proof.

b) By Proposition \ref{Proposition 2}  the condition
$\gamma_{\Omega}(y) > 0$ implies that $\gamma_{\Omega_r}(y_r) > r$.
Now we can perform calculation similar to the proof of assertion a).
\end{proof}

We need similar facts for families of sets.
Let $\{\Omega_t \} $ be a 
family of open sets in $\bR^n$.
Let $(\Omega_t)_{r}$ be $r$-neighborhood of
$\Omega_t$ in $\bR^n$. Define
\begin{equation}
\label{defEr}
E_r = \cup_t [(\Omega_t)_r \times \{ t \}], \;\;\;
{\bf \Gamma}_r = \cup_t[ \partial (\Omega_t)_r \times \{ t \}].
\end{equation}

\begin{proposition}
\label{Proposition 4t}
Let $\{\Omega_t \} $ be a locally
Lipschitz continuous 
family of open sets in $\bR^n$. 
Let for every $ t\in \bR^1_+$ the set $\Omega_t$ satisfy Condition
\ref{Cond.1}
with radius $r_0 > 0$.
Let $0<r \leq r_0$. Fix
$y\in \partial \Omega_{t_0}$. Let $y_r\in \partial
(\Omega_{t_0})_r$ be such point that $|y-y_r|=r$.  Then:
  
a) Let ${\bf \Gamma}$, defined by (\ref{defE}), be (2,1)
differentiable at $(y, t_0)$. That is
in a suitable coordinate system $(x_1, ..., x_n)$ in $\bR^n$ we have
$$ 
y=0,\; \; \Omega
\cap
B^{n+1}_\delta(0)=\{ x_n>\Psi(x',t)\} \cap
B^{n+1}_\delta(0)
$$ 
for small $\delta>0$ and 
$$
\Psi(x',t)=v(t-t_0)+
\sum^{n-1}_{i=1} \kappa_ix^2_1+o(|t-t_0|+|x'|^2).
$$
Then the point $y_r$ is unique given the point $y$,
and in the coordinate system introduced above 
$y_r=(0,-r)$,  
$$
\Omega_r\cap B^{n+1}_{\delta_1}(y_r)=\{ x_n>\Psi_r (x',t)\} 
\cap B^{n+1}_{\delta_1}(y_r) 
$$ 
for sufficiently small $\delta_1>0$, where
$$
\Psi_r(x',t)=-r+v(t-t_0)+\sum^{n-1}_{i=1}\frac{\kappa_i}
{1+\kappa_ir}x^2_1+o(|t-t_0|+|x'|^2).
$$ 

b) If ${\bf \Gamma}_r$ is (2,1)  
differentiable at $(y_r,t_0)$
and $\gamma_{\Omega_{t_0}}(y) > 0$ then  
${\bf \Gamma}$ is (2,1)  differentiable at $y$ and 
the functions $\psi$ and $\psi_r$ have same expansions at $(0, t_0)$
as in a).
\end{proposition}
\begin{proof}
The proof is similar to the one of Proposition \ref{Proposition 4}.
In particular, the sets $A_\varepsilon$, $N_\varepsilon$ and 
$\tilde{N}_{\varepsilon,\delta}$ are defined as following.
$$
A_\varepsilon=\{\;(z',z_n, t)\;\; |\;\;  z_n\geq v(t-t_0) + \sum^{n-
1}_{i=1}(\kappa_i-\varepsilon)z^2_i - \varepsilon |t-t_0| \;\}.
$$
$N_\varepsilon$ is defined as following: for each $t=t^*$ 
the set $N_\varepsilon \cap \{ t=t^* \}$ is
$\varepsilon$-neighborhood in $\bR^n$ of the set 
$\{x \;\; | \;\; (x,t^*) \in A_\varepsilon\}$. Finally, 
$$
\tilde{N}_{\varepsilon,\delta}=\{\;(z',z_n)\;\; |\;\;
z_n\geq v(t-t_0) + 
\sum^{n-1}_{i=1}(\frac{\kappa_i-\varepsilon}{1+(\kappa_i-
\varepsilon)r}-\delta)z^2_1-r -( \varepsilon + \delta) |t-t_0|\; \}.
$$ 
The rest of the argument does not change.
\end{proof}

To derive further properties of the distance function of a set $\Omega$  
satisfying the 
Condition  \ref{Cond.1},  we need to integrate over such sets. 
To do this it is convenient to decompose the set by 
suitable subsets, 
and define a local coordinate system in each subset.
The construction is carried out below and 
can be roughly described as following.
The subsets are unions of distance rays passing through
subsets of the boundary. The coordinate systems consist of 
variables $x_1, ..., x_{n-1}$ on the 
boundary, and the variable $x_n$
along the distance rays. However, since $\partial \Omega$
is not smooth enough, we have to use $\partial \Omega_r$.
Now we turn to the construction.

By Proposition \ref{Proposition 3} the set
$\Omega_r$ for small enough $r>0$ has $C^{1,1}$ boundary. Fix
such $r$. We can
choose sets ${\cal U}_1,...,{\cal U}_N$, bounded and open in $\bR^n$,
 so
that
$\partial \Omega_r \subset \cup_{k=1}^{N}{\cal U}_k$, and
 for each $k=1,...,N$ in a suitable  coordinate system
in $\bR^n$ we have
$$
\Omega_r \cap {\cal U}_k=
\{\; (x',x_n)\;\;|\;\;x'\in \tilde{U}_k,\;x_n>\Phi_k(x')\;\}\cap
{\cal U}_k,
$$ 
where $\tilde{U}_k =\{\;x'\;|\;(x',x_n)\in {\cal U}_k \;\}
 \subset \bR^{n-1}$ is an open set, and
$\Phi_k$ is a $C^{1,1}$ function on $\bR^{n-1}$.

We can also choose open sets ${\cal U}_0$ and $ {\cal U}_{N+1}\subset \bR^n$ 
such that
${\cal U}_0\subset \Omega_r$, ${\cal U}_{N+1}
\subset \bR^n\backslash \Omega_r$, both ${\cal U}_0$ and ${\cal U}_{N+1}$
do not intersect $\partial \Omega$, and
$\bR^n=\cup^{N+1}_{k=0}{\cal U}_k$.  Let $\tilde{\Psi}_k,\;k=1,...,N$, be a
smooth partition of unity on $\bR^n$ related to the sets
${\cal U} _0, {\cal U} _1,..., {\cal U} _{N+1}$, i.e., 
$\sum^{N+1}_{k=0}\tilde{\Psi}_k\equiv 1$ on
$\bR^n$ and $\tilde{\Psi}_k\in C^\infty_0({\cal U}_k)$.  Then 
by our choice of ${\cal U} _0$ and
${\cal U} _{N+1}$ we get $\sum^N_{k=1}\tilde{\Psi}_k\equiv 1$ on $\partial
\Omega_r$.

Define maps $G_k:\tilde{U}_k\times \bR^1\rightarrow
\bR^n\;\;(k=1,...,N)$ by
\begin{equation}
\label{I2}
G_k(x',x_n)=y+x_nDd^s_{\Omega_r}(y)
\end{equation}
where $y=(x',\Phi_k(x'))\in \partial \Omega_r$, and 
$d^s_{\Omega_r}(\cdot)$ is the signed distance to 
$\partial \Omega_r$. 

For each
$z\in \partial \Omega_r$ the vector $Dd^s_{\Omega_r}(z)$ is the inner
unit normal to $\partial \Omega_r$ in $z$.  Since 
$\partial \Omega_r$ is a $C^{1,1}$  manifold, 
$Dd^s_{\Omega_r}$ is Lipschitz on $\partial
\Omega_r$. Then since $\Phi_k$ is $C^{1,1}$ on $\bR^{n-1}$,
the map $G_k$ is Lipschitz on bounded subsets of
$\tilde{U}_k\times \bR^1$. In particular, the map
$\tilde{U}_k  \rightarrow \partial \Omega_r$ 
defined by $x' \rightarrow G_k(x',0)$ is Lipschitz.

From that by explicit computation we get the following
\begin{lemma}
\label{jacobian lemma}
Let the map $G_k:\tilde{U}_k\times \bR^1\rightarrow
\bR^n$ be defined by (\ref{I2}).
Then the Jacobian $JG_k$ is a locally bounded measurable function.
For ${\cal L}^{n-1}$ a.e. point $x' \in \tilde{U}_k$
the surface $\partial \Omega_r$ is twice differentiable
at the point $G_k(x',0)$. For every such point $x'$
the map $G_k$ is differentiable
at $(x', x_n)$ for every $x_n \in \bR^1$. At every
such point $(x', x_n)$ the Jacobian of $G_k$ is given 
by
\begin{equation}
\label{jacobian}
J_nG_k(x',x_n)=\sqrt{1+|D\Phi_k(x')|^2}\prod^{n-
1}_{i=1}\left[1-\kappa_{r,i} x_n\right],
\end{equation}
where $\kappa_{r,1},..., \kappa_{r,n-1}$ are the
principal curvatures of $\partial \Omega_r$ at
$G_k(x', 0)$.
\end{lemma}

Denote by $U_k$ the sets 
\begin{equation}
\label{defUk}
U_k =
\cup_{x'\in \tilde{U}_k}\{\; (x',x_n)\;|\;x_n\in
[0,\gamma_{\Omega_r}(y))\;\}\subset \bR^n,
\end{equation}
 where 
$\gamma_{\Omega_r}:\bR^n\rightarrow
\bR^1$ is the ridge function (\ref{defGamma}) of the set $\Omega_r$.
Note that  
$\gamma_{\Omega_r}(z)\geq r$ for all $z\in
\overline{\Omega}_r$ by Proposition \ref{Proposition 2} (c, d).

The map $G_k$ is one-to-one on $U_k$ (by definition of $\gamma_{\Omega_r}$
and Lemma 3.4 of \cite{EH}).  In addition   $\Omega
\backslash {\cal R}\subset \cup^N_{k=1} G_k(U_k),$  since every point of
$\Omega_r \backslash {\cal R}$ lies inside some distance ray. For each
$k=1, ..., N$ define a function
$\Psi_k:\bR^n\rightarrow \bR^1$ by
extending the function $\tilde{\Psi}_k $ from the boundary inside
$\Omega_r$ as a constant along  distance
rays and defining  $\Psi_k$ to be the $0$ on the ridge of
$\Omega_r$ and outside $\Omega_r$, i.e.,
\begin{equation}
\label{defPsi}
\Psi_k(z)=\left\{ \begin{array}{lcl}\tilde{\Psi}_k(y), &
\mbox{where }y\in {\cal N}_{\partial \Omega_r}(z), &\;\;\;
\mbox{if }z\in \Omega_r \setminus {\cal R} \\ 0, &   & \;\;\;\mbox{otherwise.}
\end{array} \right. 
\end{equation}
Note that $y$ is uniquely defuned by $z$ in the first case of (\ref{defPsi}).
We have  
$$
\mbox{supp} \Psi_k \subset G_k(U_k); \;\;\;\;
\sum^N_{k=1}\tilde{\Psi}_k (y(z))\equiv 1 \;\; \mbox{for all}\;\;
z\in \Omega_r \backslash {\cal R}. 
$$

We thus decomposed the set 
$\Omega_r \setminus {\cal R}$ by the 
sets $G_k(U_k)$ and defined a related 
partition of unity $\Psi_k$. Note that in 
general the sets $U_k$ and $G_k(U_k)$ are neither open nor closed, 
and the functions $\Psi_k$ are discontinuous. In the next propositions we
prove that these sets and functions are measurable.

\begin{proposition}
\label{Proposition I1}
Define a function
$\tilde{\gamma}_{r,k}:\bR^{n-1}\rightarrow \bR^1$ by
\begin{equation}
\label{I3}
\tilde{\gamma}_{r,k}(x')=\left\{
\begin{array}{ll}\gamma_{\Omega_r}(G_k(x',0)),
& \mbox{if }x'\in \mbox{cl}(\tilde{U}_k )\\ 0, & \mbox{if }x'\notin 
\mbox{cl}(\tilde{U}_k)
\end{array} \right.
\end{equation}
where $\mbox{cl}(\cdot)$ denotes closure of a set. Then $\tilde{\gamma}_{r,k}$
is uppersemicontinuous (usc).  The sets $U_k$ are ${\cal L}^n$ measurable.
The ridge set of $\Omega_r$ has ${\cal L}^n$ measure 0.
\end{proposition}

\begin{proof}
It  is enough to show that for any sequence $\{x_m'\}_{m=1}^{\infty}$ 
such that  $x_m'\in 
\mbox{cl}(\tilde{U}_k)$, $x_m'\rightarrow x'$
and $\tilde{\gamma}_{r,k}(x_m')\rightarrow L$ we have
$\gamma_{r,k}(x')\geq L$.

The function $\gamma_{\Omega_r}$ is usc on $\bR^n$.  
Since $x \in \mbox{cl}(\tilde{U}_k)$ and
$G_k$ is Lipschitz
on $\mbox{cl}(\tilde{U}_k) \times {0})$, 
$$
\tilde{\gamma}_{r,k}(x')=\gamma_{\Omega_r}(G_k(x', 0))
\geq \lim_{m\rightarrow \infty}
\gamma_{\Omega_r}(G_k((x_m')).
$$
Thus, $\tilde{\gamma}_{r,k}$ is usc.

We have $U_k=\{\; (x',x_n)\;|\;x'\in \tilde{U}_k,\;0\leq x_n<
\tilde{\gamma}_{r,k}(x')\; \}$. Thus $U_k$ is
${\cal{L}}^n$-measurable (note that $\tilde{U}_k$ is open in
$\bR^{n-1}$).

To prove that the ridge has ${\cal{L}}^n$ measure 0, we note that
$G_k^{-1}({\cal R})\cap
\overline{U}_k=\{ (x',x_n)|x'\in \tilde{U}_k,\;x_n=\tilde{\gamma}_{r,k }(x') \}$. 
Then ${\cal{L}}^n(G_k^{-1}({\cal R})\cap
\overline{U}_k)=0$  because 
$\tilde{\gamma}_{r,k }(\cdot)$ is
usc. Since $G_k$ is Lipschitz on 
$\mbox{cl}(\tilde{U}_k) \times [0, \diam{\Omega_r}]$
and since $\cup_{k=1}^N \overline{U}_k = \overline{\Omega}_r$, 
we have ${\cal{L}}^n({\cal R})=0$.
\end{proof}

We will use repeatedly the following fact:
\begin{lemma}
\label{Lemma I1} Let $f:A\subset \bR^m\rightarrow \bR^n (m\leq
n)$ be a Lipschitz function, $A$ be a ${\cal{L}}^m$-measurable set, 
and let
$g:\bR^m\rightarrow \bR^1$ be a ${\cal{L}}^m$-measurable function.
Then the function $g\circ f^{-1}:\bR^n \rightarrow \bR^1$ is
${\cal{L}}^n$-measurable.
\end{lemma}

\begin{proof}  We can assume that $g \equiv 0$ outside $A$.  Let
$V\subset {\bR^1}\cup \{ + \infty \}$ be an open set.
 Then $W=g^{-1}(V)\cap A$ is
${\cal{L}}^m$-measurable.  Then $(g\circ f^{-1})^{-1}V=f(g^{-
1}(U))\cap A=f(W),$ which is a measurable set since $f$ is Lipschitz
and $m\leq n$ (see \cite{EGar}, Lemma 2 of 3.3.1).
\end{proof}

\begin{proposition}
\label {Proposition I2}  
Functions $\Psi_k$ and $\Psi_k\circ G_k$
are ${\cal{L}}^n$-measurable.
\end{proposition}
\begin{proof}  The functions $G_k$ are Lipschitz. Then by 
Lemma \ref{Lemma I1} it is enough to prove that
functions $\Psi_k\circ G_k$ are ${\cal{L}}^n$-measurable  
(since then
$\Psi_k=(\Psi_k\circ G_k)\circ G^{-1}_k$).
Let $(x',x_n)\in U_k.$ Then
$$
\Psi_k(G_k(x',x_n))=
\tilde{\Psi}_k(G_k(x',0))=\tilde{\Psi}_k(x',\Phi_k(x')).
$$
Consider the function 
$\Theta(x',x_n)=\tilde{\Psi}_k(x',\Phi_k (x'))$ .
Since $\tilde{\Psi}_k\in C^\infty(\bR^{n})$ and 
$\Phi_k\in C^{1,1}(\bR^{n-1})$, 
we see that
$\Theta$ is a $C^{1,1}$
function $\bR^n\rightarrow \bR^1$. Now since  $\Psi_k\circ
G_k= \Theta$ on $U_k$,  $\Psi_k\circ
G_k \equiv 0$ outside $U_k$, and by 
Proposition \ref{Proposition I1} the set $U_k$ 
 is measurable,
the assertion follows.
\end{proof}

\begin{proposition}
\label{Lemma I2}  
Let $\Omega$ be a bounded open set satisfying the Condition \ref{Cond.1}. 
Let $X\subset
\partial \Omega$ and ${\cal H}^{n-1}(X)=0$.  Let $Y$ be the
union of distance rays passing 
through $X$, i.e., $Y=\cup_{x\in X}R_x$.
Then ${\cal L}^n(Y)=0$.
\end{proposition}
\begin{proof}
We can assume that $X\subset G_k(U_k)\cap \partial \Omega$
for some $k$, and
then $Y\subset G_k(U_k)$. We also can assume
that
\begin{equation}
\label{canAssume}
\gamma_{\Omega}(x) > 0 \;\;\; \mbox{ for all } \; x \in X.
\end{equation}
If (\ref{canAssume}) is not true then we can replace $X$ by
$X \cap \{ x \; | \; \gamma_{\Omega}(x) > 0 \}$ and $Y$ does
not change.

 Since $\tilde{\gamma}_{r,k}(x')$ is a
u.s.c. function of $x'$, we get ${\cal L}^n(\{ x_n=\gamma_{r,k}(x',0)\})
=0$.  We also know that $J_nG_k>0$ on the set $U_k$.
The map $G_k$ is one-to-one on $U_k$, and
so $N(G_k,y)=1$ for $y\in Y \subset  G_k(U_k)$, 
where $N(f,z)$ denotes the multiplicity function, which is the number of 
elements of $f^{-1}(z)$.  Then by area
formula we obtain
\begin{equation}
\label{arF1}
{\cal L}^n(Y) = \int_{\bR^n}\chi _Y(y) dy=
\int_{U_k}\chi _Y(G_k(x))J_nG_k(x)dx,
\end{equation}
where $\chi_A(\cdot)$ is the characteristic function of the set
$A$, i.e., $\chi_A(\cdot)$ equals 1 on $A$ and 0 outside $A$.
By Lemma \ref {jacobian lemma} the Jacobian  $J_nG_k$
is a locally bounded measurable function.
Let $x' \in \tilde{U}_k$ be such point that
the surface $\partial \Omega_r$ is twice differentiable
at the point $G_k(x',0)$.
By Lemma \ref {jacobian lemma} the Jacobian  $J_nG_k$
is 
defined by  the expression (\ref{jacobian}) at  
$(x', x_n)$ for every $x_n \in \bR^1$. The
principal curvatures of $\partial \Omega_r$ at the point
$G_k(x', 0)$ satisfy
\begin{equation}
\label{kappa-gamma}
 \kappa_{r,i}>-\frac{1}{r}, \;\;\;
\kappa_{r,i}\tilde{\gamma}_{r,k}(x')\leq 1
\end{equation}
by definition of $\tilde{\gamma}_{r,k}(\cdot)$,
 and 
by Condition  \ref{Cond.1}. We also have
\begin{equation}
\label{gammaBd}
\tilde{\gamma}_{r,k}(x') \leq \mbox{diam}\Omega_r.
\end{equation}
Let $m \leq n$ be such that
$-\frac{1}{r} \leq \kappa_{r,1},..., \kappa_{r,m} < 0$ and 
$\kappa_{r,m+1},..., \kappa_{r, n-1} \geq 0$.
Then for $x_n \in (r,\tilde{\gamma}_{r,k}(x')) $ 
we calculate using  (\ref{kappa-gamma}), (\ref{gammaBd}):
\begin{equation}
\label{jacIneq}
0 \leq \frac{J_nG_k(x', x_n)} {J_nG_k(x', r)} \leq
\frac{\prod^{m}_{i=1}\left[1-\kappa_{r,i} x_n\right]}
{\prod^{m}_{i=1}\left[1-\kappa_{r,i} r\right]} \leq
\left[1+\frac{1}{r}\tilde{\gamma}_{r,k}(x')\right]^m
\leq C,
\end{equation}
where $C$ depends on $r$, $n$ and $\mbox{diam}(\Omega)$.
By Lemma  \ref {jacobian lemma} the inequalities (\ref{jacIneq})
hold for a.e. $x \in U_k \cap \{(x',x_n) \;|\;x_n>r\}$.

Define a map
$\tilde{G}_k:\tilde{U}_k\rightarrow 
\partial \Omega \subset \bR^n$ by
$\tilde{G}_k(x')= G_k(x',r) $. Then by (\ref{jacobian})
\begin{equation}
\label{JtoJ}
J_{n-1}\tilde{G}_k (x')=J_nG_k(x',r)
\;\; \mbox{for a.e.} \;x'\in\tilde{U}_k.
\end{equation}
So using the fact that the relation $Y\subset \overline{\Omega}$
 implies $\chi_Y(G_k(x',x_n))=0$ for
$x_n<r$, we get from (\ref{arF1}), (\ref{gammaBd}), (\ref{jacIneq}) and 
(\ref{JtoJ}) using area formula:
\begin{eqnarray}
{\cal L}^n(Y) &\leq& C\int_{U_k}\chi_Y(G_k(x))\;J_{n-
1}\tilde{G}_k(x')dx
\nonumber \\
 & = & C\int_{\tilde{U}_k}J_{n-
1}\tilde{G}_k(x')[\int_r^{\tilde{\gamma}_{r,k}(x')}
\chi_Y(G_k(x',x_n))dx_n]dx' 
\nonumber \\
 & = & C \int_{\tilde{U}_k} \chi_X(\tilde{G}_k(x')) 
[\tilde{\gamma}_{r,k}(x')-
r]\cdot J_{n-1}\tilde{G}_k(x')dx'
\nonumber \\
 & \leq &
C\int_{\tilde{U}_k}\chi_X(\tilde{G}_k(x'))\;J_{n-
1}\tilde{G}(x')dx'
\nonumber \\
 & = &
C\int_{\partial \Omega}\chi_X(z)N(\tilde{G}_k,z) d{\cal H}^{n-
1}(z) = C{\cal H}^{n-1}(X)=0,
\nonumber
\end{eqnarray}
where we have used the fact that $N(\tilde{G}_k,z)= 1$
for $z\in X$ (which follows from (\ref{canAssume}) and Proposition 
\ref{Proposition 2} (b,c) ).
\end{proof}
\begin{remark}
Without assuming a condition of the type of Condition \ref{Cond.1}
the assertion of Proposition \ref{Lemma I2}
is not true. For example consider a nonconvex polygon in
$R^2$ and take $X$ to be one point, a vertex of a re-enterant corner.
\end{remark}
\begin{remark}
\label{InvMeasZero}
Note that it follows from the proof above that
$$
{\cal L}^n(G_k^{-1}(Y)) = 
\int_{\tilde{U}_k}J_{n-
1}\tilde{G}_k(x')[\int_r^{\tilde{\gamma}_{r,k}(x')}
\chi_Y(G_k(x',x_n))dx_n]dx' = 0.
$$
\end{remark}

The following proposition describes the structure of the boundary of a
set that satisfies Condition \ref{Cond.1}.
\begin{proposition}
\label{Proposition 4.1}
Let $\Omega$ be a bounded set.  Let $\Omega$ 
satisfy the condition \ref{Cond.1} with radius $r_0$. Then 
$\partial \Omega$ is $({\cal H}^{n-1}, n-1)$ rectifiable subset of
$\bR^n$. In addition, $\partial \Omega$ is twice differentiable
${\cal H}^{n-1}$ a.e. on the set 
\begin{equation}
\label{set B}
{\cal B}=\{ x \in \partial \Omega \; | \; \gamma_{\Omega}(x)>0 \}.
\end{equation}
\end{proposition}
\begin{proof}
Let $r>0$ be such number that $\partial \Omega_r$ is $C^{1,1}$.
Existence of such $r$ follows from Proposition \ref{Proposition 3}.
It follows from Proposition \ref{Proposition 2} b)
that the nearest point projection 
mapping $P: \Omega_r \rightarrow \Omega$ 
is well-defined and onto. Moreover, by 
inequality 4.8(8) of \cite{FedCMeas}
(or by inequality (\ref{1.1a}) proved
in Appendix \ref{Apndx 1} below)
the map $P$ is Lipschitz. Since $\partial \Omega_r$ is
$C^{1,1}$, we get that $\partial \Omega$ is  
$({\cal H}^{n-1}, n-1)$ rectifiable.

By Proposition \ref{Proposition 2} every point of the set ${\cal B}$
has a unique nearest point on $\partial \Omega_r$. Since 
$\partial \Omega_r$ is
$C^{1,1}$, it follows that
$\partial \Omega_r$ is twice differentiable ${\cal H}^{n-1}$ a.e.
Since the map $P$ is Lipschitz it follows that
 for  ${\cal H}^{n-1}$  a.e. 
$y \in {\cal B}$  
the surface $\partial \Omega_r$ is twice differentiable
at the corresponding 
point $y_r$. Applying Proposition \ref{Proposition 4} b), 
we conclude the proof.
\end{proof}

Propositions \ref{Proposition I1}, \ref{Lemma I2} and \ref{Proposition 4.1}
imply 

\begin{corollary}
\label{A.E.}
Let $\Omega \subset \bR^n$ satisfy Condition \ref{Cond.1}.
Then ${\cal L}^n$ a.e point of $\Omega$ lies in the relative
 interior of
a distance ray that intersects  $\partial \Omega$ at
 a point at
which  $\partial \Omega$ is twice differentiable.
\end {corollary}

We need an analogue of the Proposition \ref{Proposition 4.1} for
a Lipschitz family $\{\Omega_t\}$ of sets that satisfy  
Condition \ref{Cond.1} with a radius $r_0$ independent of $t$.
Let $E, \;
{\bf \Gamma} $ be defined by (\ref{defE}).
Let $E_r, \;
{\bf \Gamma}_r $ be defined by (\ref{defEr}).
Denote
\begin{equation}
\label{defFxt}
D(x,t) = d^{s}_{\Omega_t} (x),
\end{equation}
where $d^{s}_{\Omega_t} (x)$ is the signed distance 
function (\ref{sgnDist}).

\begin{proposition}
\label{Proposition 4.1.t}
Let $\{\Omega_t\}$ be a locally Lipschitz family of open 
bounded sets that 
satisfy the Condition \ref{Cond.1} with radius $r_0$ independent
of $t$. Then  

a) ${\cal L}^{n+1}$ a.e.
$(x_0, t_0) \in E$
 lies in the relative interior of
a distance ray of the set $\Omega_{t_0}$
that intersects $\partial \Omega_{t_0}$
 at such point $y \in \partial \Omega_{t_0}$ 
that ${\bf \Gamma}$ is (2,1) differentiable at $(y,t_0)$;

b) There exists a constant $C$ depending only on $r_0$
such that the function $D(x,t) - C|x|^2$ 
on the set $E_{r_0}$ is concave in 
$x$ for every $t$.

\end{proposition}
\begin{proof}
First we prove the following fact.

\begin{lemma}
\label{semiconcDist}
Let $\Omega$ satisfies Condition  \ref{Cond.1} with radius $r_0$.
Then the distance function $d_{\Omega}(x)$ is semiconcave
in $\Omega$,
i.e. there exists a constant $C$ depending only on $r_0$ such
that the function $d_{\Omega}(x) - C|x|^2$ is concave in $\Omega$.
\end{lemma}
\begin{proof}
Let $x\in \Omega$. Let $w \in \bR^n$ satisfy
$$
|w| \leq d_{\Omega}(x).
$$ 
Let $y \in {\cal N}_{\partial \Omega_{r_0}}(x)$, where 
$\Omega_{r_0}$ is $r_0$- neighborhood of
$\Omega$. Then by Proposition \ref{Proposition 2}(c),
\begin{eqnarray}
|x-y|& =& d_{\Omega}(x) + r_0, 
\nonumber \\
|x+w-y| & \geq & d_{\Omega}(x+w) + r_0,
\nonumber \\
|x-w-y| & \geq & d_{\Omega}(x-w) + r_0.
\nonumber
\end{eqnarray}
Then we have
$$
d_{\Omega}(x+w) - 2d_{\Omega}(x) + d_{\Omega}(x-w)  \leq
|x+w-y| - 2|x-y| + |x-w-y| \leq \frac{C}{r_0}|w|^2,
$$
and so 
$d_{\Omega} -2 \frac{C}{r_0}|x|^2$ is concave.
\end{proof}

Now we prove Proposition \ref{Proposition 4.1.t}.
 The function $D(x,t)$ defined in (\ref{defFxt}) is
Lipschitz since $\{\Omega_t\}$ is a Lipschitz 
family of sets. By Proposition \ref{Proposition 2} the 
set $(\Omega_t)_{r_0}$
satisfies Condition  \ref{Cond.1} with radius $\frac{r_0}{2}$,
and 
$$
d^{s}_{(\Omega_t)_{r_0}} (x) = d^{s}_{\Omega_t} (x) + r_0
\;\;\;\mbox{for }\;(x,t) \in E_{r_0}. 
$$
Then by Lemma \ref{semiconcDist} the function  $D(x,t) - C|x|^2$ 
on the set $E_{r_0}$ is concave in 
$x$ for every $t$.
Now by  Theorem 1 of
Appendix 2  of \cite{Krylov} the function $D(x,t)$ is (2,1) 
differentiable
at ${\cal L}^{n+1}$ a.e. point of $E_{r_0}$. It follows
from Proposition \ref{Proposition 4t}
that if $D(\cdot)$ is (2,1) differentiable at a point
$(x_0, t_0) \in E$, then ${\bf \Gamma}$ is (2,1) differentiable 
at $(y,t_0)$, where 
$y \in {\cal N}_{\partial \Omega_{t_0}}(x_0)$. 
Note that $x_0$ lies in the relative interior
of its distance ray if $d_{\Omega_{t_0}}$ is twice $x$-differentiable
at $x_0$. The Proposition is proved.
\end{proof}

\section{Properties of mass transport density.}
\label{SectDens}
Let $\Omega$ be a bounded set that satisfies the condition \ref{Cond.1}
with radius $r_0$.
Let $t>0$.
Define a function $a=a_{\Omega, t}:\bR^n\rightarrow \bR^1$ as following.  Let
$x\in \partial \Omega$ be such a point that
  the surface $\partial \Omega$ is twice differentiable at 
$x$ and let
$\kappa_1, ..., \kappa_{n-1}$ be the principal curvatures of  $\partial  
\Omega$ at
$x$.  Introduce the
coordinate $s$ on the distance ray $R_x$ 
by $s(z)=|z-x|$ for $z\in R_x$. Then
$s$ changes in the interval $(0, \gamma(x))$.
Define the function $a$ on $R_x$ as the following:
\begin{equation}
\label{4}
a(s)={1\over t}\prod^{n-1}_{i=1}{1\over {1-{s
\kappa_i}}}\int^s_0 \prod^{n-1}_{i=1}(1-\xi \kappa_i)\left[
\frac{\int^{\gamma(x)}_0 \nu \prod^{n-1}_{i=1} (1-\kappa_i
\nu)d\nu}{\int^{\gamma(x)}_0 \prod^{n-1}_{i=1} (1-\kappa_i
\nu)d\nu}-\xi \right] d\xi
\end{equation}

By  Corollary  \ref{A.E.} the function $a$ is now 
defined a.e. in $\Omega$.

Define $a\equiv 0$ on
$\bR^n \backslash \Omega$.
Now $a$ is defined a.e. in $\bR^n$.

\begin{definition}
\label{MasTrDen}
The function $a(\cdot)=a_{\Omega, t}(\cdot)$ is called 
mass transport density.
\end{definition}

Define the function $V: \partial \Omega \rightarrow \bR$
at the point of twice differentiability of $\partial \Omega$
as following.
\begin{equation}
\label{Veloc}
V(x)=
{1\over
t}\frac{\int^{\gamma (x)}_0 s\prod^{n-1}_{i=1} (1-\kappa_i
s)ds}{\int^{\gamma (x)}_0 \prod^{n-1}_{i=1} (1-\kappa_i
s)ds}.
\end{equation}

\begin{proposition}
\label{Proposition 1}
Let $\Omega \in \bR^n$ be a bounded open set that satisfies 
Condition \ref{Cond.1} with radius $r_0$.
Then there exists a constant $C$ depending only on 
$r_0,\; n,\; \diam\,{\Omega}$ such that the following is true.
Let $x \in \partial  \Omega$ be such point that 
$\partial \Omega$
be twice differentiable at $x$.
Then the mass transport density on the distance ray
$R_x$, defined by (\ref{4}) for
 $s\in [0,\gamma(x)]$, satisfies: \\
(a) $|a(s)|\leq C$. 
Derivative $a'(s)$ exists and is continuous for $s \in
(0,\gamma(x))$. \\
(b) $a(s)$ satisfies the equation 
$$
a'(s)-a(s)\sum^{n-1}_{i=1} 
{{\kappa_i}\over
{1-\kappa_i s}} -{s \over t}+V(x)=0
$$ 
on
$(0,\gamma(x))$,  the boundary conditions
$a(0)=0$, $a(\gamma(x))=0$,
and $a(s)>0$ on $(0,\gamma(x)).$
\\
(c) For $s \in (0,\gamma (x))$ the inequalities hold
$$
|a(s)|\leq C(\gamma(x)-
s),\;\;\;\;
|a(s)|\leq Cs.
$$ 
(d) $|a'(s)|\leq C$ on $(0,\gamma(x)).$

\end{proposition}

\begin{proof}
(a) Existence and continuity of $a'(s)$ on $(0,\gamma(x))$ are
checked explicitly.  The bound $|a(s)|\leq C$ will follow from
$(c)$. \\

(b) The equation is checked explicitly.  The boundary conditions
follow from $(c)$. The inequality  $a(s)>0$ on $(0,\gamma(x))$
holds since 
either $\kappa_i < 0$, or $0\leq
\gamma(x)\leq {1\over \kappa_i}$ and then, defining
$$
\tilde{V}(s)=\frac{\int^s_0 \nu\prod^{n-1}_{i=1}(1-\kappa_i
\nu)d\nu}{\int^s_0 \prod^{n-1}_{i=1}(1-\kappa_i \nu)d\nu}
$$
we calculate:
$$
{d\tilde{V} \over ds}=\frac{s\prod^{n-1}_{i=1}(1-
\kappa_i s) \int^s_0 \prod^{n-1}_{i=1}(1-\kappa_i \nu)d\nu-
\prod^{n-1}_{i=1}(1-\kappa_i s)\int^s_0 \nu\prod^{n-1}_{i=1}(1-
\kappa_i \nu)d\nu}{[\int^s_0 \prod^{n-1}_{i=1}(1-\kappa_i \nu)d\nu)]^2}
$$
$$
=\frac{\prod^{n-1}_{i=1}(1-\kappa_i s)}{[\int^s_0 \prod^{n-1}_{i=1}(1-
\kappa_i \nu)d\nu)]^2}\int^s_0 \prod^{n-1}_{i=1}(1-\kappa_i  
\nu)(s-\nu)d\nu>0\
\ \mbox{for } s\in (0,\gamma(x)),
$$
and so
$$
a(s)={1\over t}
[\tilde{V}(\gamma(x))-\tilde{V}(s)]
\prod^{n-1}_{i=1}{1\over {1-\kappa_i
s}}\int^s_0 \prod^{n-1}_{i=1}(1-\kappa_i \nu)d\nu
>0\ \ \mbox{for }s \in
(0,\gamma(x)).
$$

(c)We write $\gamma$ for $\gamma(x)$:
Denote $\Phi(s)=\int^s_0 \nu\prod^{n-1}_{i=1}(1-
\kappa_i\nu)d\nu,\;\Psi(s)=\int^s_0
\prod^{n-1}_{i=1}(1-\kappa_i\nu)d\nu$. For $s \in (0,\gamma)$  
we calculate:
\begin{eqnarray}
ta(s) & = &
\frac{[\Phi(\gamma)\Psi(s)-\Phi(\gamma)\Psi(\gamma)]+
[\Phi(\gamma)\Psi(\gamma)-\Phi(s)\Psi(\gamma)]}{\Psi (\gamma)}
\prod^{n-1}_{i=1}{1\over {1-\kappa_is}}
\nonumber \\
& = & (\gamma-s)
\frac{-\Phi(\gamma)\Psi'(s_1)+\Phi'(s_2)\Psi(\gamma)}{\Psi
(\gamma)}\prod^{n-1}_{i=1}\frac{1}{1-\kappa_i
s},
\nonumber
\end{eqnarray}
where $s_1,s_2 \in [s,\gamma]$. Note, that
$\tilde{V}(\xi)={\Phi(\xi)\over {\Psi(\xi)}}$.  So
\begin{eqnarray}
ta(s)& = & \prod^{n-1}_{i=1}{1\over {1-\kappa_is}}\left[-
\tilde{V}(\gamma) \prod^{n-1}_{i=1}(1-\kappa_is_1)+\prod^{n-
1}_{i=1}(a-\kappa_is_2)\cdot s_2\right]\;(\gamma-s)
\nonumber \\
& = &
\left[-\tilde{V}(\gamma)\prod^{n-1}_{i=1}\frac{1-\kappa_is_1}{1-
\kappa_is}+s_2\prod^{n-1}_{i=1}\frac{1-\kappa_is_2}{1-
\kappa_is}\right]\;(\gamma-s).
\nonumber
\end{eqnarray}
If $\kappa_i\geq 0$ then for $s^*\in(s,\gamma)$ we have $0\leq s
\kappa_i<s^*\kappa_i<1,$ and so
$$
0\leq \frac{1-\kappa_i s^*}{1-\kappa_i
s}\leq1.
$$
If $\kappa_i< 0$ then it follows from the Condition \ref{Cond.1} that
$|\kappa_i| \leq C$, and so 
$$
0 \leq \frac{1-\kappa_i s^*}{1-\kappa_is} = 
\frac{1+|\kappa_i| s^*}{1+|\kappa_i|s} \leq 1+C\gamma 
\leq 1+ C \diam\Omega.
$$
It remains to estimate $\tilde{V}(s)$ for $s \in [0, \gamma]$.
\begin{equation}
\label{tildVEst}
\tilde{V}(s)\leq
\frac{\int_0^s \gamma
\prod^{n-1}_{i=1} (1-\nu\kappa_i)d\nu}{\int_0^s \prod^{n-
1}_{i=1} (1-\nu\kappa_i)d\nu}=\gamma \leq \mbox{diam}\Omega.
\end{equation}
So, $ta(s)\leq C\mbox{diam}\;\Omega\;(\gamma-s),$ and, since 
$t>0$,
\begin{equation}
\label{5}
0<a(s)\leq C(\gamma-s)\mbox{ for }s \in (0,\gamma(x)).
\end{equation}
To estimate $a(s)$ near $s=0$, we compute using (\ref{tildVEst})
\begin{eqnarray}
|ta(s)| & = & 
\left|[\tilde{V}(\gamma)- \tilde{V}(s)]
\prod^{n-1}_{i=1}\frac{1}{1-\kappa_i s}\int^s_0
\prod^{n-1}_{i=1}(1-\kappa_i\nu)d\nu \right|
\nonumber \\
& \leq & \gamma \left|\prod^{n-1}_{i=1}\frac{1}{1-\kappa_i
s}\int^s_0 \prod^{n-1}_{i=1}(1-\kappa_i\nu)d\nu\right|=\gamma
s \prod^{n-1}_{i=1}\frac{1-\kappa_i s^*}{1-\kappa_i s},
\nonumber
\end{eqnarray}
where $s^*\in[0,s]$.  
Since
$\kappa_i \gamma \leq 1$, then for
$s \leq {\gamma \over 2}$ we get 
$$
0 \leq \frac{1-\kappa_i s^*}{1-\kappa_i
s}\leq 2 \;\;\mbox{ if } \kappa_i \geq 0,
$$
and since 
$\kappa_i > -\frac{1}{r}$ we get 
$$
0 \leq \frac{1-\kappa_i s^*}{1-\kappa_i
s}\leq C(\mbox{diam}(\Omega), r) \;\;\mbox{ if } \kappa_i < 0.
$$
So
 $|ta(s)|\leq C \gamma s$,
and, since $ t>0$, we get 
$$
|a(s)|\leq C s.
$$
If $s\geq{\gamma \over 2}$, then $s\geq\gamma-s$, and from
(\ref{5}) 
$$
|a(s)|\leq C(\gamma-s)\leq Cs.
$$
\\

(d) We prove that $|a'(s)|\leq C$ for $s \in (0,\gamma)$.
$$
a'(s)=a(s)\sum^{n-1}_{i=1}\frac{\kappa_i}{1-\kappa_i
s}+{s \over t}-V(x).
$$
  By (\ref{tildVEst}), $|V(x)|\leq
\frac{\gamma}{t}$.  
It remains to prove
that 
\begin{equation}
\label{MTrDenDerivEst}
\left|a(s)\sum^{n-1}_{i=1}\frac{\kappa_i}{1-\kappa_i
s}\right|\leq C 
\end{equation}
If $\kappa_i \geq 0$ then $0\leq \kappa_i \gamma\leq 1$,
and since $\gamma>s$ we get:
$$
0\leq\frac{\kappa_i}{1-\kappa_i s}\leq\frac{\kappa_i}{\kappa_i
\gamma-\kappa_i s}=\frac{1}{\gamma-s},
$$
If $\kappa_i < 0$, then 
$$
\left|\frac{\kappa_i }{1-\kappa_i s}\right|\leq \frac{1}{s}.
$$
Now (\ref{MTrDenDerivEst}) follows from (c).
\end{proof}

\section{Mass balance equation}
\label{ChngVar}

Let $\{\Omega_t\}$ be a Lipschitz family of 
open sets satisfying Condition
\ref{Cond.1} with radius $r_0>0$. Let $\{\Omega_t\}$ 
satisfy equation (\ref{1}) on the 
time-interval $[a,b]$ in the sense described
in Theorem \ref{Theorem 3}.
Let the function $w(x,t)$ be defined by (\ref{2}).
Let $a(x,t)$ be the function $a_{\Omega_t, t}(x)$  
from Definition \ref{MasTrDen}.

The purpose of this section is to show that for a.e. $t \in [a,b]$
\begin{equation}
\label{I1}
\int_{\bR^n} aDwD\phi dz=\int_{\bR^n} ({w\over t}-w_t)\phi
dz
\end{equation}
for any $\phi \in C^\infty(\bR^n)$, where $w_t$ denotes $\partial_t w$.
Thus we show that the function $a(\cdot, t)$ satisfies
(\ref{defMassTrDen}) for a.e $t$.

It follows from Proposition \ref{Proposition 4.1.t} 
that for a. e. $t$ the following is true:
${\cal L}^n$ a. e. point $x\in \Omega_t$ lies in
the relative interior of the distance ray that 
intersects $\partial \Omega_t$ at
such point $y\in \partial \Omega_t$ that the surface ${\bf \Gamma}$
is (2,1) differentiable at $(y,t)$, where ${\bf \Gamma}$ is defined
by (\ref{defE}).  Fix such $t$.
For the rest of this section
 we drop $t$ in the notation, i.e., we 
write $\Omega$, $a(x)$, $w(x)$, $w_t(x)$ for 
$\Omega_{t}$, $a(x, t)$, $w(x, t)$, $w_t(x, t)$.

We use same notation $a(\cdot)$ for both mass transport density
$a(x)$, a function defined on  $\bR^n$, and for mass transport density
on a ray $R_x$, the  function  $a(s)$ defined on $[0, \gamma(x)]$.

Denote by $\tilde{\Gamma}$ the subset of $\partial \Omega$ that 
consists of all points at which the surface
${\bf \Gamma}$ is $(2,1)$-differentiable.
By choice of $t$ and Proposition \ref{Proposition 4.1} we get
$$
{\cal H}^{n-1}({\cal B}  \setminus \tilde{\Gamma})=0,
$$
where the set ${\cal B}$ is defined by (\ref{set B}).
From definition of ${\cal B}$ it follows that no distance
rays of $\Omega$ have their lower ends in the subset
$\partial \Omega \setminus {\cal B}$ of the boundary.
Denote by $R(\tilde{\Gamma})$  the union of all distance rays that have
lower ends in  $\tilde{\Gamma}$, i.e., 
$R(\tilde{\Gamma})= \cup_{x \in \tilde{\Gamma}}R_x $. Then it follows from 
Proposition \ref{Lemma I2}
and 
Remark \ref{InvMeasZero}
that
\begin{equation}
\label{MeasZeroRel}
{\cal L}^{n} ( \Omega \setminus R(\tilde{\Gamma}) ) = 0, \;\;\;
{\cal L}^{n}[ G_k^{-1} (\Omega \setminus R(\tilde{\Gamma}))] =0 \; 
\mbox{for}\;k=1,...,N.
\end{equation}

Let $V(\cdot)$ be the function
(\ref{Veloc}). It follows from (\ref{MeasZeroRel}) 
that at a. e. point $x$ of $\Omega$ the function $w_t(x)$ is given
by the expression 
$\frac{1}{t}V(y)$ where 
$y \in {\cal N}_{\partial \Omega}(x)$. 
In the next proposition we prove that $w_t$  is measurable.

We continue to use local coordinate
systems on $\Omega_r$, described in Section \ref{Sets}
and notation introduced there, in particular
sets $U_k$, partition of unity
$\Psi_k$  and coordinate mappings $G_k$ defined 
by (\ref{defUk}), (\ref{defPsi}) and (\ref{I2}) respectively.

\begin{proposition}
\label{Proposition I3}

a) For $k=1,...,n$ there exists a bounded 
${\cal{L}}^n$ measurable
function $v_k:U_k\rightarrow \bR^1$ such that
\begin{equation}
\label{I5}
v_k(x',x_n)=V(G_k(x',r))\;\;\; \mbox{if } \;\;\;
G_k(x',r)\in \tilde{\Gamma},\; x_n \in [r,\tilde{\gamma}_{r,k}(x')].
\end{equation}

b) The function $w_t$ is ${\cal L}^n$ measurable and for any
bounded
measurable function $\varphi:\; U_k\rightarrow \bR^1$
\begin{equation}
\label{I6}
{1\over t}\int_{U_k\cap \{x_n\geq r \}}
v_k{\varphi}J_nG_kdx=\int_{G_k(U_k)}w_t(\varphi \circ G^{-1}_k)dx
\end{equation}
\end{proposition}
\begin{proof}

a) $\partial \Omega_r$ is a $C^{1,1}$ manifold, 
and so the second fundamental 
form of  $\partial \Omega_r$ is defined as a differential form with
$L^{\infty}$ coefficients. At all points where  $\partial \Omega_r$
is twice differentiable the principal curvatures of 
$\partial \Omega_r$ are the roots of the characteristic polynomial
of the second fundamental form.  To write this in coordinates, 
we recall that
$\partial \Omega_r \cap G_k(U_k)$ is the graph of 
a $C^{1,1}$ function $x_n = \Phi_k(x')$ defined 
on the set $\tilde{U}_k \in \bR^{n-1}$.
Then there exist 
$L^{\infty}$  functions $[u_{i j}(x')]_{i,j=1}^n$
such that at the points where $D^2\Phi_k$ exists,
$\frac{\partial^2 \Phi_k}{\partial x_i \partial x_j}
= u_{ij}$.

Denote by $P_{x'}(s)$ the characteristic polynomial of the matrix
$[u_{ij}(x')]$ i.e. 
$P_{x'}(s)=\mbox{det}([u_{ij}(x')]-sI_{n-1})$.
Note that the coefficients of the polynomial $P_{x'}(\cdot)$
are measurable functions of  $x'$ since they are compositions of 
polynomials and measurable functions $u_{ij}(x')$.

The
eigenvalues of the matrix $[u_{ij}(x')]_{i,j=1}^n$
are principal curvatures of $\partial \Omega_r$
at the points of twice differentiability of $\Phi_k$.
Let $x' $ be such a point,
let $\kappa_{r,1},...,\kappa_{r,n-1}$ be principal curvatures of
$\partial \Omega_r$ in $G_k(x',0)$, and let $m$ be such number
 that $\kappa_{r,i}\neq 0$ for $i=1,...,m$ and $\kappa_{r,i}=0$
for $i=m+1,...,n-1$.
Then $P_{x'}(s)=\prod^{n-1}_{i=1}(\kappa_{r,i}-s)$, and so
$$
(-1)^{n-1}s^{n-1}P_{x'}({1\over s})=\prod^{n-1}_{i=1}\;\;(1-
\kappa_{r,i}s)=
\prod^m_{i=1}\kappa_{r,i} \prod^m_{i=1}
(\frac{1}{\kappa_{r,i} }-s).
$$

For any $x'\in U_k$ the point $G_k(x',r)$ belongs to
$\partial \Omega$ and $|G_k(x',0)-G_k(x',r)|=r$.
If the surface $\partial \Omega$ is twice differentiable at the
point $y=G_k(x',r)$, then, by 
Proposition \ref{Proposition 4},
the surface $\partial \Omega_r$ is twice differentiable at the
point $y_r=G_k(x',0)$,
and, denoting by $\{ {\kappa}_i \}_{i=1}^n$ and
 $\{ {\kappa}_{r,i} \}_{i=1}^n $ the principal
curvatures of $\partial \Omega$ at $y$ and of $\partial \Omega_r$
at $y_r$, we get
$$
{\kappa}_{r,i} =\frac{{\kappa}_i}{1+{\kappa}_ir},\;\mbox{so
}\frac{1}{{\kappa}_{r,i} }=\frac{1}{{\kappa}_i}+r\;\;(\mbox{if  
}{\kappa}_i\neq 0).
$$
Then, assuming  as before that  
${\kappa}_1,...,{\kappa}_m \neq 0$ and  
${\kappa}_{m+1},...,{\kappa}_{n-1}=0$, we get
if  $m \ne 0$
\begin{eqnarray}
\prod^{n-1}_{i=1}(1-{\kappa}_is)
&=& \prod^m_{i=1}{\kappa}_i\prod
^m_{i=1}(\frac{1}{{\kappa}_i}-s)
=\prod^m_{i=1}{\kappa}_i\prod
^m_{i=1}\left[\frac{1}{{\kappa}_{r,i}}-(s+r)\right]
\nonumber \\
 & = &
(-1)^{n-
1}\frac{\prod^m_{i=1}{\kappa}_i}{\prod^m_{i=
1}{\kappa}_{r,i}}(s+r)^{n-1}P_{x'}(\frac{1}{s+r}).
\nonumber
\end{eqnarray}
If $m=0$ (i.e., all $\kappa_{r,i}\equiv 0$), then
$\prod^{n-1}_{i=1}(1-{\kappa}_is)=(-1)^{n-1}(s+r)^{n-
1}P_{x'}(\frac{1}{s+r})=1$.  

Now we see that at all $x'$ as above,
the function $V$ defined in  (\ref{Veloc}) can be expressed  at 
the point $ G_k(x',r) \in \partial \Omega$ as
\begin{equation}
\label{I7}
V(G_k(x',r))=\frac{\int^{{\tilde{\gamma}_{r,k}}(x')-r
}_0
s(s+r)^{n-
1}P_{x'}(\frac{1}{s+r})ds}{\int^{{\tilde{\gamma}_{r,k}}(x')-r}_0
(s+r)^{n-1}P_{x'}(\frac{1}{s+r})ds}
\end{equation}
since, by Proposition \ref{Proposition 2},  
$\gamma_{\Omega_r}(z)=\gamma(z)_{\Omega}+r$ for $z\in
\Omega$.

Define a function $\eta=\eta(x',x_n)$ on  $\tilde{U_k} \times \bR^1 $ by the 
right-hand side of (\ref{I7}).
This function $\eta$ is ${\cal L}^n$-measurable:  
$\eta$ does not depend on $x_n$, the coefficients
of the polynomial $P_{x'}(\cdot)$ are measurable
functions of $x'$ and by 
Proposition \ref{Proposition I1},
the functions $\tilde{\gamma}_{r,k}$ are measurable
functions of  $x'$.

Then it follows from Proposition \ref{Proposition I1}
 that the 
function 
$$v_k(x',x_n)=\eta (x',x_n) \chi_{U_k}(x',x_n)$$
satisfies all properties asserted in (a) (note that
the right-hand side of (\ref{I5}) is bounded by diam $\Omega$).

By Lemma 
\ref{Lemma I1} the
function $v_k\circ G^{-1}_k$ is ${\cal L}^n$-measurable.  
By (\ref{MeasZeroRel}) we get 
\begin{equation}
\label{I8'}
\begin{array}{l}
{1\over t}v_k\circ G^{-1}_k=w_t\;\;\;\mbox{at}\; {\cal L}^n 
\;\mbox{a.e. point of}\;\,
G_k(U_k \cap \{ \, x_n \geq r \,\}), 
 \\
 \\
\mbox{and } \;\; w_t \equiv 0 \;\;\; on \;\;\; 
G_k(U_k \cap \{ \, x_n < r \,\})
\end{array}
\end{equation}
and so $w_t$ is measurable.  

Since map $G_k$ is one-to-one and Lipschitz on
$U_k$,
and since $U_k$ is measurable,
we get by
Theorem 3.2.5 of \cite{Federer} and Lemma \ref{Lemma I1}
that for every bounded measurable function $\varphi$ 
\begin{equation}
\label{I8}
\int_{U_k}v_k\varphi \;J_nG_kdz=\int_{G_k(U_k)}(v_k\circ
G^{-1}_k)(x)
(\varphi \circ G^{-1}_k)(x) dx
\end{equation}

Now (\ref{I6}) follows from (\ref{I8}), (\ref{I8'}).
\end{proof}

\begin{corollary}
\label{Corollary I1}  
The function ${w\over t}-w_t$ is
measurable, and for each $\varphi \in C^\infty(\bR^n)$
\begin{equation}
\label{I9}
\int_{\bR^n}({w\over t}-w_t)\varphi dz={1\over
t}\sum^N_{k=1}\int_{G_k^{-1}(R(\tilde{\Gamma}) ) }
(w\circ G_k-v_k)(\varphi \circ
G_k)(\Psi_k\circ G_k)J_nG_kdx.
\end{equation}
\end{corollary}

\begin{proof}
We have
 $${w\over t}-w_t=\sum^N_{k=1}({w\over t}-
w_t)\Psi_k.
$$
  Since $\mbox{supp}\;\Psi_k\subset G_k(U_k)$, the integrand
on the
right-hand side of (\ref{I9}) is measurable by 
Propositions \ref{Proposition I3}  and
\ref{Lemma I2}.  It is also clear that the integrand on the
right-hand side
is bounded.

Then from Proposition \ref{Proposition I3}  (b) 
and area formula, we get
$$
\int_{\bR^n}({w\over t}-w_t)\varphi dz={1\over
t}\sum^N_{k=1}\int_{\bR^n}(w\circ G_k-v_k)(\varphi \circ
G_k)(\Psi_k \circ G_k)J_nG_k dx.
$$
This formula is true for each $\varphi \in C^\infty$, and, by
approximation, for every bounded measurable $\varphi$.  Now
(\ref{MeasZeroRel})
implies (\ref{I9}).
\end{proof}

Now we prove that the mass transport
density $a(\cdot)$ is measurable. The proof is
similar to the proof of Proposition \ref{Proposition I3}.
\begin{proposition}
\label{Proposition I4}
The function  $a(\cdot)$ is  ${\cal L}^n$-measurable function
on $\bR^n$. The function  
$a \circ G_k^{-1}$ is  ${\cal L}^n$-measurable function
on $U_k$.
\end{proposition}

\begin{proof}
For $(x', s)$ such that 
$$
G_k(x',r)\in \tilde{\Gamma},\;\;\;s\in
[r,\tilde{\gamma}_{r,k}(x',r)],
$$
 the function $a\circ G_k(x',s)$
is given  by  expression (\ref{4}), where the curvatures are 
computed at the point  $G_k(x',r)$. Calculations
similar to those in Proposition \ref{Proposition I3} show
that  this expression
can be rewritten using the characteristic polynomial
$P_{x'}(\cdot)$ of the matrix $D^2\Phi_k(x')$ as
following:
$$
 a\circ G_k(x',s)={1\over t}\;\frac{\int^{s-r}_0
(\nu+r)^{n-
1}P_{x'}(\frac{1}{\nu+r})d\nu}{(s+r)^{n-1}P_{x'}(\frac{1}{s+r})}\times
$$
\begin{equation}
\label{I12}
\left[ \frac{\int^{\tilde{\gamma}_{k,r}(x')-r}_0 \nu(\nu+r)^{n-
1}P_{x'}(\frac{1}{\nu+r})d\nu}{\int^{\tilde{\gamma}_{k,r}(x')-r}_0
(\nu+r)^{n-1}P_{x'}(\frac{1}{\nu+r})d\nu}-\frac{\int^{s-r}_0 \nu(\nu+r)^{n-
1}P_{x'}(\frac{1}{\nu+r})d\nu}{\int^{s-r}_0 (\nu+r)^{n-
1}P_{x'}(\frac{1}{\nu+r})d\nu} \right]
\end{equation}
By (\ref{MeasZeroRel}), the function  
$a \circ G_k^{-1}(\cdot)$ is defined by the above formula
at a.e. point of $U_k$, and is
${\cal L}^n$-measurable function of
variables $x',s$ on $U_k$  since  it is a rational function of variables
$s$ and $\tilde{\gamma}_{k,r}(x')$
 with measurable coefficients 
and $\tilde{\gamma}_{k,r}(x')$
is ${\cal L}^{n-1}$ measurable function of $x'$.
By Lemma 
\ref{Lemma I1}, the function $a(\cdot)$ is measurable.
\end{proof}

{\bf Proof of (\ref{I1}).}
We  will transform the right-hand side of (\ref{I9}). 
From (\ref{defPsi}) we get
\begin{equation}
\label{PsiIndepOfs}
\Psi_k(G_k(x',s)) = \Psi_k(G_k(x',r)), \;\;\;\;
s \in (r, \tilde{\gamma}_{r,k}(x')).
\end{equation}
Let $\tilde{G}_k: \tilde{U}_k \rightarrow \partial\Omega$
be the map defined by
$$
\tilde{G}_k(x') = G_k(x',r).
$$
 Then for
$k=1,...,N,$ using definitions of $w,v_k,\Psi_k$ we compute
\begin{eqnarray}
{1\over t}\int_{G_k^{-1}(R(\tilde{\Gamma}) )}
(w\circ G_k 
- v_k)(\varphi \circ
G_k)(\Psi_k \circ G_k)J_nG_k dx 
\hspace{8em} & &
\nonumber \\
 = \int_{\tilde{G}_k^{-1}(\tilde{\Gamma}) }
\Psi_k(G_k(x',r)) \left\{ 
\int^{\tilde{\gamma}_{r, k}(x')}_r {1\over t}[s-r 
\begin{array}{l}  \\ \\ \end{array} \right. \hspace{3em} & &
\label{RHS} \\
 & & \nonumber \\
  \hspace{5em}
\left. \begin{array}{l} \\ \\ \end{array}
-v_k(x',s)]\varphi (G_k(x',s))J_nG_k(x',s)ds
 \right\} dx'  & &
\nonumber
\end{eqnarray}

Let $x' \in \bR^{n-1}$ be such that $G_k(x', r) \in \tilde{\Gamma} $.  
Then $\partial \Omega$ is
twice differentiable at $G_k(x',r)$, and so 
$$
{1\over t}v_k(x',s)= V(G_k(x',r))\;\;\;
\mbox{ for }\;s\in [r, \tilde{\gamma}_{r,k}(x')].
$$ 
By Proposition  \ref{Proposition 4}, $\partial \Omega_r$ is twice 
differentiable at
$G_k(x',0)$.  Let $\kappa_1,...,\kappa_{n-1}$ be principal
curvatures of $\partial \Omega$ at $G_k(x',r)$, and
$\kappa_{r,1},...,\kappa_{r,n-1}$ be principal curvatures of
$\partial \Omega_r$ at $G_k(x',0)$. Then
$\kappa_{r,i} =\frac{\kappa_i}{1+r\kappa_i}$. Let
$y=G_k(x',r)$. Then  the mass transport density
 on the
ray $R_y$ satisfies ODE from Proposition \ref{Proposition 1}(b).
Thus we get for $s\in
[r,\tilde{\gamma}_{r,k}(x')]=[r,\gamma_{\Omega}(y)+r]$:
\begin{eqnarray}
{1\over t}(s-r-v_k)
& = &
-\frac{da(s-r)}{ds}+a(s-r)\sum^{n-
1}_{i=1}\frac{\kappa_i}{1-\kappa_i(s-r)}
\nonumber \\
& = &
-\frac{da(s-r)}{ds}+a(s-
r)\sum^{n-1}_{i=1}\frac{\kappa_{r,i} }{1-\kappa_{r,i} s}
\nonumber
\end{eqnarray}
where $a(s)$ is mass transport density on the
ray $R_y$.

By Proposition \ref{Proposition 1}, $\frac{da}{ds}$ exists and is
continuous and  
bounded on the interval
$(0,\tilde{\gamma}_{r,k}(x')-r)$, and
$$
a(0)=a(\tilde{\gamma}_{r,k}(x')-r)=0.
$$
We also see by
(\ref{jacobian})  that
\begin{equation}
\label{jacDeriv}
\frac{\partial J_nG_k}{\partial s}(x',s)=J_nG_k(x',s)\sum^{n-
1}_{i=1}\frac{-\kappa_{r,i} }{1-\kappa_{r,i} s}.
\end{equation}
Then we calculate, integrating by parts:
$$
\int^{\tilde{\gamma}_{r, k}(x')}_r {1\over t}[s-r-
v_k(x',s)]\varphi (G_k(x',s))J_nG_k(x',s)ds 
 = 
$$
$$
\int^{\tilde{\gamma}_{r,k}(x')}_r\left[-\frac{da(s-r )}{ds}-a(s-
r)\frac{\frac{\partial}{\partial
s}J_nG_k(x',s)}{J_nG_k(x',s)}\right]\varphi (G_k(x',s))J_nG_k(x',s)ds
$$
$$
=\int^{\tilde{\gamma}_{r,k}(x')}_r a(s-
r)\frac{\partial}{\partial
s}\varphi(G_k(x',s))J_nG_k(x',s)ds
$$
This allows to derive from (\ref{I9}), (\ref{RHS}) and
(\ref{PsiIndepOfs}) the following equality:
\begin{eqnarray}
\label{I10}
\int_{\bR^n}({w\over t}-w_t)\varphi dz   =
\hspace{22em}
& &
\\
\sum_{k=1}^N 
\int_{\tilde{G}_k^{-1}(\tilde{\Gamma})}
\int^{\tilde{\gamma}_{r,k}(x')}_r 
\Psi_k(G_k(x',s)) 
a(G_k(x',s))\frac{\partial}{\partial s}\varphi
(G_k(x',s)J_nG_k(x',s)ds dx' &&
\nonumber
\end{eqnarray}

Let $z=G_k(x',s)$. Using the fact that 
distance function is twice differentiable on $R(\tilde{\Gamma})$, we get:
\begin{eqnarray}
\frac{\partial}{\partial s}
(\varphi(G_k(x',s)) &=& ( D\varphi)(G_k(x',s))\cdot 
\frac{\partial }{\partial s}(G_k(x',s))
\nonumber \\
&=& D\varphi (z)\cdot Dd_\Omega(z).
\label{s-deriv}
\end{eqnarray}

We insert the right-hand side of (\ref{s-deriv})
 into the right-hand side of (\ref{I10}) and 
use area formula
and (\ref{MeasZeroRel})  to change
variables from $(x', s)$ to $z$ in the integral on the
right-hand side of (\ref{I10}). Then using that 
$\sum_k
\Psi_k\equiv 1$ on $\Omega^0 \backslash {\cal R}$ and that
$a \equiv 0$ on ${\cal R}$ and in $\bR^n \setminus \Omega$,
we get (\ref{I1}).

\section{Proof of Theorem \ref{Theorem 3}}
\label{Th.3}
To conclude the proof  of Theorem \ref{Theorem 3} we
need to demonstrate the following.
Let, as above, $w(x)=d_{\Omega}(x)$, and let 
$v\in W^{1,\infty}(\bR^n)$ satisfy $|Dv|\leq 1$ a.e. Then
\begin{equation}
\label{I15}
\int_{\bR^n}({w\over t}-w_t)w\geq \int_{\bR^n}({w\over t}-w_t)v.
\end{equation}

We show that this follows from  (\ref{I1}).
By approximation, (\ref{I1}) is true for Lipschitz $\varphi$.
Let $\varphi=w-v$.  Then by (\ref{I1}) that it is
enough to prove that
$$\int_{\bR^n} aDw(Dw-Dv)\geq 0.$$
But, since $a\geq 0$ in $\bR^n$ 
 and $a \equiv 0$ outside $\Omega^0$ and on ${\cal R}$, and
$|Dw|\equiv 1$ on $\Omega \backslash {\cal R}$,
$$
\int_{\bR^n} aDw Dw=\int_{\bR^n}a\geq \int_{\bR^n}aDw
Dv
$$
since  $ |DwDv|\leq\sqrt{|Dw|^2|Dv|^2}\leq 1$.

Theorem \ref{Theorem 3} is proved.

\section{Proof of Theorem \ref{Theorem 4}}
\label{Th.4}

The proof of Theorem \ref{Theorem 4} 
is similar to the proof  Theorem \ref{Theorem 3}. 
We give a formal calculation. 
Each step can be justified the same way as
in the proof of Theorem \ref{Theorem 3}.

Fix $t\in [a, b]$. Define
$$
{\cal D}_1={\cal D}_1(t)=\{ x\in
\Omega^1_t\; |\; d_{\Omega^1_t}(x)>d_{\Omega^2_t}(x)\},
$$
$$
{\cal D}_2={\cal D}_2(t)=\{ x\in
\Omega^2_t \; |\; d_{\Omega^2_t}(x)>d_{\Omega^1_t}(x)\}.
$$
Let $r>0$, let $\Omega^1_r=\Omega^1_{t,r}$ and
$\Omega^2_r=\Omega^2_{t,r}$ be $r$-neighborhoods of
$\Omega^1=\Omega^1_t$ and $\Omega^2=\Omega^2_t$,
and let $r$ be so small that $\partial \Omega^1_r$ and 
$\partial \Omega^2_r$ are $C^{1,1}$ manifolds.

Let l=1,2. Let $\{ \tilde{U}^l_k\}^{N_l}_{k=1}$ and 
$\{ U^l_k\}^{N_l}_{k=1}$ be the
coordinate neighborhoods associated the set $\Omega^l_r$
 same way as $\{ \tilde{U}_k\}$, $\{ U_k\}$ with $\Omega_r$ in the
proof of Theorem \ref{Theorem 3}. Let $G^l_k:U^l_k\rightarrow
\bR^n$ be the corresponding coordinate
mappings (\ref{I2}) and let $\{ \Psi^l_k\}^{N_l+1}_{k=0}$
 be the corresponding partitions of unity (\ref{defPsi}).
Let $\tilde{\gamma}^l_{r,k}$ be the functions defined 
as in (\ref{I3})  from the functions 
$\gamma_{\Omega_r^l}$ 
using the coordinate mappings  $G^l_k$ for $k=1,...,N_l$.
Define the following functions 
$\tilde{\delta}^l_{r,k}$  for $k=1,...,N_l$. Let
$\delta (\cdot)\equiv 
 \gamma_{\Omega^1  \cap \Omega^2}(\cdot)$.
Then
$$
\tilde{\delta}^l_{r,k}(x')=\left\{
\begin{array}{ll}\delta(G^l_k(x',r))+r,
& \mbox{if }x'\in \mbox{cl}(\tilde{U}^l_k )
\;\; \mbox{ and }\;\; G^l_k(x',r) \in 
\partial(\Omega^1  \cap \Omega^2);
\\ 
r, & \mbox{otherwise}.
\end{array} \right.
$$
Note that
$$
(G^l_k)^{-1}({\cal D}_l) \cap  U^l_k = \{\; (x', s) \;\;\; |
\;\;\; x' \in \tilde{U}^l_k, \;s \in [\,\tilde{\delta}^l_{r,k}(x'),\;
\tilde{\gamma}^l_{r,k}(x')\,] \;\}.
$$

The functions $\tilde{\delta}^l_{r,k}$ are uppersemicontinuous,
the proof is similar to the proof of Proposition \ref{Proposition I1},
using the properties of sets satisfying Condition  \ref{Cond.1}.

Define the mass transport density $a(\cdot)$.
Let  $a\equiv 0$ on $\bR^n \backslash (\Omega^1 \cup \Omega^2)$.
To define $a$ in $\Omega^1 \cup \Omega^2$ it is enough to define
$a$ on ${\cal D}_1 \cup {\cal D}_2$.
Let $x\in {\cal D}_l\backslash {\cal R}_l$, where $l=1,2$.
Let $y \in {\cal N}_{\partial \Omega_l}(x)$
and $s = |x-y|$. Then we define
$$
a(x)={1\over t}\prod^{n-1}_{i=1}\frac{1}{1-s
\kappa_i}\int^s_{\delta(y)} \prod^{n-1}_{j=1}(1-\kappa_j\xi)\left[
\frac{\int^{\gamma^l(y)}_{\delta(y)}\nu \prod^{n-1}_{k=1}(1-
\kappa_k\nu)d\nu}{\int^{\gamma^l(y)}_{\delta(y)} \prod^{n-1}_{k=1}(1-
\kappa_k\nu)d\nu}-\xi \right]d\xi,
$$
where $\kappa_k$ are curvatures of $\partial \Omega^l_t$ at $y$.
Then $a$ satisfies the equation (\ref{TrDenODE}) on
subintervals of distance rays of ${\Omega^l_t}$ that lie 
in ${\cal D}_l,$ and
$a\equiv 0$ on $\partial {\cal D}_l$ and on
${\cal R}_l,\;l=1,2$.  We calculate, denoting by $\kappa_{r,i} $ the
principal
curvatures of $\partial \Omega^l_r$ at the  point $y_r \in \Omega^l_r$ 
nearest to $x$ 
\begin{eqnarray}
\int_{\bR^n}({w\over t}&-&w_t)\varphi dz
\nonumber \\
&=&{1\over
t}\sum^2_{l=1}\sum^{N_l}_{k=1}\int_{U^l_k\cap (G^l_k)^{-1}({\cal
D}_l)}(w\circ G^l_k-
V\circ G^l_k)(\varphi \circ G^l_k)(\Psi^l_k
\circ G^l_k)J_nG^l_k\,dx
\nonumber \\
&=&\sum^2_{l=1}\sum^{N_l}_{k=1}
\int_{\tilde{U}^l_k}
\Psi^l_k(G^l_k(x',0))
\nonumber \\
& & \times 
\int^{\tilde{\gamma}^l_{r,k}(x')}_{\tilde{\delta}^l_{r,k}(x')}
{1\over t}\left[(s-r)-
V(G^l_k(x',0))\right](\varphi(G^l_k(x',s))J_nG^l_k(x',s)\, ds dx'
\nonumber \\
&=&\sum^2_{l=1}\sum^{N_l}_{k=1}
\int_{\tilde{U}^l_k}
\Psi^l_k(G^l_k(x',0)) 
\int^{\tilde{\gamma}^l_{r,k}(x')}_{\tilde{\delta}^l_{r,k}(x')}\left( -
\frac{d\;a(G^l_k(x',s))}{d\;s} \right.
\nonumber \\
& & \left.
+a(G^l_k(x',s))\sum^{n-
1}_{i=1}\frac{\kappa_{r,i} }{1-\kappa_{r,i} s} \right)
\varphi \, J_nG^l_k(x',s)\,ds dx'
\nonumber \\
&=& \sum^2_{l=1}\sum^{N_l}_{k=1}
\int_{\tilde{U}^l_k}
\Psi^l_k(G^l_k(x',0))
\nonumber \\
& & \times 
\int^{\tilde{\gamma}^l_{r,k}(x')}_{\tilde{\delta}^l_{r,k}(x')}
a(G^l_k(x',s))\
\frac{\partial \varphi(G^l_k(x',s ))}{\partial s}JG^l_k(x',s)\,ds dx'
\nonumber \\
&=&\sum^2_{l=1}\sum^{N_l}_{k=1}\int_{G^l_k(U^l_k)\cap {\cal
D}_l}a(z)\Psi_k^l(z)\;D\varphi (z)\;Dw(z)\,dz
\nonumber \\
&=& \int_{\bR^n}aDwD\varphi\, dz.
\nonumber
\end{eqnarray}
Thus we have proved that the mass balance equation is satisfied.
This implies Theorem \ref{Theorem 4}.

\section{Compression molding model}
\label{CompMod} 

Compression molding is the process of deformation of an incompressible
plastic material between
two horizontal plates. The following
 simplified mathematical model of the process  was derived by
G. Aronsson \cite{ArCompMol} based on Hele-Show model for 
non-Newtonian fluid. 
Suppose that the
distance between the horizontal
plates is small.
 Then we can assume that the region 
occupied by plastic at each time $t$
has the form $U_t = \Omega_t \times [0, h_t]$
where $\Omega_t \subset \bR^2$, and that the
pressure in plastic 
does not depend on the vertical coordinate,
i.e., pressure is the function $u(x,t)$
where $x \in \Omega_t$. Evolution of rescaled
$\Omega_t$ and $u(x,t)$ is described by
the following free boundary problem.
Given
an open set $\Omega_0 \in \bR^2$ find an expanding family of open
sets $\Omega_t \in \bR^2, t \geq 0,$ and a  function $u(x,t)$ defined on 
$\cup_t (\Omega_t\times \{ t \})$ such that
\begin{eqnarray}
- \mbox{div} (|Du|^{ p-2 } Du)=1 & \;\;\; & \mbox{in } \Omega_t,
\label{CompMoldEqP}  \\
u=0 &  & \mbox{on } \Gamma_t,
\label{bdCondP} \\
V= |Du|^{p-2} & &  \mbox{on } \Gamma_t, 
\label{fBdCondP}
\end{eqnarray}
where $\Gamma_t = \partial \Omega_t$, $V$ denotes the outer normal 
velocity of $\Gamma_t$.
Condition 
(\ref{fBdCondP}) means that the free boundary $\Gamma_t$ moves 
with the velocity of the flow.

In the paper \cite{AEComprMot} the asymptotic limit as $p \rightarrow
\infty$ in the problem (\ref{CompMoldEqP}) -  (\ref{fBdCondP}) 
  was considered.  This limit 
corresponds to the case of highly non-Newtonian
fluid. 

It was shown in \cite{AEComprMot} that formally sending to a limit
in (\ref{CompMoldEqP}) -  (\ref{fBdCondP}) one obtains the following 
problem. Find a family $\{ \Omega_t \}$ of open subsets  of $ \bR^2$,  a 
function $u(x,t)$,  and a mass transport density function $a(x,t)$ 
satisfying:
\begin{equation}
\label{DensLimFormOm}
\left\{
\begin{array}{lll}
|Du| \leq 1 & \;\;\; & \mbox{a.e. in }  \Omega_t,  \\
a \geq 0 & & \mbox{a.e. in }  \Omega_t, \\
\mbox{supp}(a) \in \{ |Du| = 1 \}, & & \\
- \mbox{div}(aDu)=1 & & \mbox{ in }  \Omega_t,
\end{array}
\right.
\end{equation}
where the last equation is understood in the weak sense, and
\begin{equation}
\label{DensLimFormGam}
u=0, \;\; V=a \;\; \mbox{ on } \Gamma_t, 
\end{equation}
where $V$ is the outer normal velocity of $\Gamma_t$.

It was shown in \cite{AEComprMot} that solutions of
(\ref{DensLimFormOm}), (\ref{DensLimFormGam}) have the
form
\begin{equation}
\label{DistFormula}
u(x,t)=d_{ \Omega_t}(x),
\end{equation}
where the right-hand side is defined by 
(\ref{distFunc00}), and that
formally the problem (\ref{DensLimFormOm}) - 
(\ref{DensLimFormGam}) can be rewritten as
following.
Find $\{\Omega_t \}$  such that
\begin{equation}
\label{CompMoldEv}
w- \partial_t w \in \partial I_{\infty} [u],
\end{equation}
where
\begin{equation}
\label{indFnct}
w(\cdot, t) = \chi_{\Omega_t}(\cdot),
\end{equation}
where $ \chi_{\Omega_t}(\cdot)$ is the indicator function 
of $\Omega_t$ that equals 1 inside and 0 outside $\Omega_t$.
Existence of
a weak solution of  (\ref{DistFormula}) - (\ref{CompMoldEv})
was proved in \cite{AEComprMot}.
Namely, there exists a family $\{\Omega_t\}$ of sets of finite 
perimeter such that $\partial_t w$ is a nonnegative Radon measure
and
\begin{eqnarray}
\int_{\bR^2} w(x, t)(v(x)-u(x, t))dx \leq 
\int_{\bR^2} (v-u(\cdot, t)) d(\partial_t w(\cdot,t)) 
 \nonumber \\
 \label{weakCompMold}  \\
\;\; \mbox{for every $v$ with } 
|Dv| \leq 1, \;\mbox{for a.e. } t.
 \nonumber
\end{eqnarray}
The following law of motion of the free boundary $\Gamma_t$
was derived in \cite{AEComprMot} by a formal calculation:
\begin{equation}
\label{MotBdCompMol}
V=\gamma (1-\frac{\kappa\gamma}{2}),
\end{equation}
where
$$
\left\{  
\begin{array}{l}
V = \mbox{ outer normal velocity of } \Gamma_t, \\
\gamma = \mbox{ function  $\gamma_{\Omega_t}(\cdot)$
defined by (\ref{defGamma})}, \\
\kappa =  \mbox{ curvature of } \Gamma_t.
\end{array}
\right.
$$

The equation (\ref{MotBdCompMol}) was derived as following.
Starting from  (\ref{DensLimFormOm}), we perform calculations 
similar to the ones that lead from
(\ref{defMassTrDen}) to (\ref{TrDenODE}).
Thus we deduce that equation
 (\ref{DensLimFormOm}) can be formally
rewritten on each distance ray $R_x$, where
$x \in \Gamma_t$, as ODE
\begin{equation}
\label{TrDenODECMold}
a'(s)-a(s)
\frac{\kappa}{1-\kappa s} + 1 = 0, \;\;\;
s\in(0,\gamma_{\Omega_t}(x)).
\end{equation}
By the nature of mass transfer process in the
compression molding model
(i.e., mass transfer from within
the set onto the boundary),
we expect that the mass transport density equals to zero 
on the ridge of $\Omega_t$ and equals to the outer normal
velocity at the boundary. This translates into
the following boundary conditions for
the ODE (\ref{TrDenODECMold})
on $R_x$:
\begin{equation}
\label{TrDenODECMoldBCond}
a(0) = V(x), \;\;\; a(\gamma(x))=0.
\end{equation}
The function $a(\cdot)$ and the constant $V$ can be found
from (\ref{TrDenODECMold}), (\ref{TrDenODECMoldBCond}). 
$V$ has the expression (\ref{MotBdCompMol}).

In this section we prove the connection between the
variational equation (\ref{CompMoldEv}) and the
geometric equation (\ref{MotBdCompMol}).

Let $E, {\bf \Gamma}$  be defined by (\ref{defE}).

\begin{theorem}
\label{Theorem 1-CM}
  Let $\{\Omega_t \},\; t\in \bR^1_+$, be an expanding, locally
Lipschitz continuous 
family of open bounded sets. 
Let for every $ t\in \bR^1_+$ the set $\Omega_t$ satisfy Condition
\ref{Cond.1}
with radius $r_0 > 0$.
 Suppose that the equation (\ref{MotBdCompMol}) 
is satisfied at every point  of $(2,1)$ differentiability of the surface
${\bf \Gamma}$. 
Then the equation (\ref{CompMoldEv}) is satisfied for a.e. $t\in \bR^1$,
where the functions $u(x,t), \; w(x,t)$ are defined by
(\ref{DistFormula}), (\ref{indFnct}).
\end{theorem}
\begin{proof}
The proof goes along the lines of the proof of 
Theorem \ref{Theorem 3}. We will sketch the proof and
present some details for the parts that are different from the proof of 
Theorem \ref{Theorem 3}.

The main steps of the proof are following: \\
{\it Step 1}. Definition and properties of mass transport density.\\
{\it Step 2}. Show that the main mass balance equation is satisfied. \\
{\it Step 3}. Show that the assertion of the Theorem follows from
the main mass balance equation.

We discuss each step.

{\it Step 1}.

   Let $(x, t)$ be
a point  of $(2,1)$ differentiability of the surface
${\bf \Gamma}$ in $\bR^2 \times \bR^1$. 
Define mass transport density $a(y,t)= a(s)$ on
the distance ray $R_x$ of $\Omega_t$
as the solution of
(\ref{TrDenODECMold}), (\ref{TrDenODECMoldBCond}).
Do this for every such point $(x, t)$.
It follows from the hypothesis of the Theorem and
from  Proposition \ref{Proposition 4.1.t}
that $a(y,t)$ is now defined at ${\cal L}^{3}$ a.e. point of
the set $E = \cup_{t}(\Omega_t \times \{t\} )
\subset \bR^2 \times \bR^1_+$. Define $a(y,t) $ by zero
at all other points of $\bR^2 \times \bR^1_+$.

Solving (\ref{TrDenODECMold}), (\ref{TrDenODECMoldBCond})
(and taking into account (\ref{MotBdCompMol}))
we see that on the ray $R_x$
\begin{equation}
\label{trDenRayCM}
a(s)=\frac{1-\kappa s}{2\kappa} - 
\frac{(1-\kappa\gamma)^2}{2\kappa(1-\kappa s)}.
\end{equation}
This function satisfies assertions (a) and (d) of 
Proposition \ref{Proposition 1}. Instead of 
assertion (c) of 
Proposition \ref{Proposition 1} we have
\begin{equation}
\label{lipTrDen}
|a(s)|\leq C(\gamma(x)-
s),\;\;\;\;
|a(s)- V(x)|\leq Cs.
\end{equation}
To see this, we
note that from the conditions of the Theorem and 
(\ref{MotBdCompMol}) it follows that
$$
V(x) \leq C(r_0, \diam \Omega_t).
$$
Then the assertions (a) and (d) of 
Proposition \ref{Proposition 1}
for the function $a(s)$ defined by (\ref{trDenRayCM})
follow
from (\ref{lipTrDen}) like in the proof of 
Proposition \ref{Proposition 1}. 

Let us prove (\ref{lipTrDen}). We can rewrite (\ref{trDenRayCM})
as
$$
a(s) = \frac{\gamma -s}{2}\left(1 +
\frac{1-\kappa \gamma}{1-\kappa s}\right).
$$
Now, since $0 \leq s \leq \gamma \leq \diam \Omega_t$ and
$\kappa \gamma \leq 1$ and $\kappa \geq -\frac{1}{r_0}$,
we get the estimate
\begin{equation}
\label{est1}
0 \leq \frac{1-\kappa\gamma}{1-\kappa s} \leq C(r_0,  \diam \Omega_t),
\end{equation}
and the first inequality of (\ref{lipTrDen}) follows. To prove the 
second, we calculate using (\ref{TrDenODECMoldBCond})
$$
a(s) - V = a(s)-a(0) = \frac{s}{2}\left[-1-
\frac{(1-\kappa\gamma)^2}{1-\kappa s}\right],
$$
and use (\ref{est1}) and the inequality
$$
0 \leq 1-\kappa \gamma \leq C(r_0,  \diam \Omega_t)
$$
to finish the proof of (\ref{lipTrDen}).

{\it Step 2}.

The purpose of this step is to prove that for any
smooth function $\varphi \in C^1_c(\bR^2 \times (0,T))$ we have
\begin{equation}
\label{balEqComprMold}
\int_0^T \int_{\bR^2} w (\varphi + \varphi_t)
-aDuD\varphi = 0.
\end{equation}
This is main mass balance equation
for compression molding model.

The proof follows Section \ref{ChngVar}.

Fix $t$.
We use local coordinate
systems on $(\Omega_t)_{r}$ where $0<r\leq r_0$,
defined in Section \ref{Sets}
and notation introduced there, in particular
sets $\tilde{U}_k$, $U_k$, partition of unity
$\Psi_k$, coordinate mappings $G_k$ 
and functions $\tilde{\gamma}_{r,k}$
defined 
by (\ref{defUk}), (\ref{defPsi}), (\ref{I2}) and
(\ref{I3}) respectively
(k=1,...,N).

We use notation $E_r$, ${\bf \Gamma}_r$ introduced in (\ref{defEr}).
From Proposition \ref{Proposition 2} it follows that for
$r \in [0, r_0]$
$$
{\bf \Gamma}_r = \{( x ,t)\in \bR^2 \times (0,T) \;\;
| \;\; d^s_{\Omega_t}(x) = -r \}.
$$
The function
$$
(x,t) \rightarrow d^s_{\Omega_t}(x)
$$
is Lipschitz since $\{ \Omega_t \}$ is a Lipschitz family of sets.
Then by Proposition 3.2.15 of \cite{Federer}, for a.e. $r \in [0, r_0]$
\begin{equation}
 set \;\; {\bf \Gamma}_r\;\; is\;\; ({\cal H}^2, 2)\;\; rectifiable.
\label{assm1}
\end{equation}
In addition, by Lemma \ref{semiconcDist}, the
function  $(x,t) \rightarrow d^s_{\Omega_t}(x) - C|x|^2$, 
where $C$ is large enough depending on $r_0$,
is concave in $x$ for every $t$ in the
set $E_{r_0}$.
Then it follows from Theorem 1 of
Appendix 2  of \cite{Krylov} and Propositions
\ref {Proposition 4.1.t}(b)
and
\ref{Proposition 3} 
that for a.e. $r \in [0, r_0]$
\begin{eqnarray}
 & &surface\;\;{\bf \Gamma}_r\;\; is \;\; (2,1)\;\, 
differentiable \;\, at \;\;{\cal H}^2\;\;
a.e. \;\, point;
\label{assm2}  \\
&& \nonumber \\
 & &functions \;\; \partial_t d^s_{\Omega_t}(x), \;
Dd^s_{\Omega_t}(x),\; D^2 d^s_{\Omega_t}(x)\;\, are\; \, 
 bounded\;\, and\;\;
\nonumber \\
& & {\cal H}^2\;
 measurable\; on\; {\bf \Gamma}_r. 
\label{assm3}
\end{eqnarray}

Function $(x,t) \rightarrow \gamma_{\Omega_t}(x)$
is uppersemicontinuous since the family of sets
$\{\Omega_t\}$ is Lipschitz.
Then similar to Propositions \ref{Proposition I3},
\ref{Proposition I4}  and Corollary \ref{Corollary I1}
we prove that the mass transport density 
$a(y,t)$ defined at Step 1
is Lebesgue measurable in $\bR^2 \times \bR^1_+$. In the proof
we use (\ref{assm3}) (i.e., we choose such $r$ in the
definition of local coordinate systems that (\ref{assm3}) is
satisfied).

It follows from Proposition \ref{Proposition 4.1.t} 
that for a.e. $t$ we have the following:
a.e. point $x \in \Omega_t$ lies in the  relative interior of
the distance ray of the set $\Omega_t$ that intersects
$\partial \Omega_t$ at a
point of (2,1) differentiability of the surface 
${\bf \Gamma}$.
Fix such $t$. Let $r\in (0, r_0]$.
Let $y \in \partial \Omega_t$ be a point
of (2,1) differentiability of ${\bf \Gamma}$. Let $\kappa^t$
be the curvature of $\partial \Omega_t$ at $y$,
let  denote $\kappa_r^t$ be the curvature of $\partial( \Omega_t)_r$
at the unique  point $y\in \partial (\Omega_t)_r$ nearest to $x$,
then
$$
\frac{1}{\kappa^t_r} = \frac{1}{\kappa^t} + r.
$$
We calculate using properties 
of the function $a$ on the rays proved in the Step 1,
equation  (\ref{TrDenODECMold}),
and Lemma \ref{jacobian lemma} with $n=2$:
\begin{eqnarray}
\int_{\bR^2}w \varphi \,dz &=&\int_{\Omega_t} \varphi \,dz 
\nonumber \\
 &=& 
\sum^{N}_{k=1}\int_{U_k}(\varphi \circ G_k)(\Psi_k
\circ G_k)J_2G_kdx
\nonumber \\
 &=& \sum^{N}_{k=1}
\int_{\tilde{U}_k}
\Psi_k(G_k(x',r))
\int^{\tilde{\gamma}_{r,k}(x')}_r
\left[\frac{da(s-r)}{ds} \right.
\nonumber \\
 & & \left.  +a(s-r)\frac{\kappa^t_{r} }{1-\kappa^t_{r} s}\right]
(\varphi(G_k(x',s))J_2G_k(x',s) ds dx'.
\nonumber 
\end{eqnarray}
In the last expression we integrate by parts and use
(\ref{jacDeriv}), (\ref{TrDenODECMoldBCond}) and 
(\ref{s-deriv}). Then we get

\begin{eqnarray}
\int_{\bR^2}w \varphi \,dz &=&
\sum^{N}_{k=1}
\int_{\tilde{U}_k}
\Psi_k(G_k(x',r)) 
\left[ \begin{array}{l}  \\ \\ \end{array} \right.
a(0)\varphi(G_k(x',r) )J_2G_k(x',r)
\nonumber \\
& & 
+  \int^{\tilde{\gamma}_{r,k}(x')}_r
a(s-r)
\frac{\partial \varphi(G_k(x',s ))}{\partial s}J_2G_k(x',s)ds
\left. \begin{array}{l}  \\ \\ \end{array} \right ]
dx'.
\nonumber \\
&=&   
\int_{\hat{\partial}( \Omega_t)_r }\varphi(P^t_r(y)) 
 \partial_t d^s_{\Omega_t}(P^t_r(y)) (1-r \kappa_r(y)) d{\cal H}^1(y)
+ \int_{\Omega_t}aD \varphi Du \,dz,
\nonumber
\end{eqnarray}
where
\begin{equation}
\label{gammaBdry}
\hat{\partial}( \Omega_t)_r = \{ y\in\partial ( \Omega_t)_r \;\; |
\gamma_{( \Omega_t)_r}(y) > r \},
\end{equation}
and 
\begin{equation}
\label{nearPtPr}
P^t_r : (\Omega_t)_r \rightarrow \Omega_t
\end{equation}
is the nearest point projection 
mapping. By Proposition \ref{Proposition 2} b) the mapping $P^t_r$
is well-defined and onto.

By Proposition \ref{Proposition 4t}, if ${\bf \Gamma}_r$ is 
(2,1) differentiable
at $y \in \hat{\partial}( \Omega_t)_r$, then ${\bf \Gamma}$ is
(2,1) differentiable
at $y'=P^t_r(y)$ and 
$$
\partial_t d^s_{\Omega_t}(P^t_r(y)) = \partial_t d^s_{\Omega_t}(y)
$$
By the equality 
$$
d^{s}_{(\Omega_t)_{r}} (x) = d^{s}_{\Omega_t} (x) + r
\;\;\;\mbox{for }\;(x,t) \in E_{r_0}, \;\; r \in (0, r_0]
$$
we see that 
$$
\partial_t d^s_{\Omega_t}(P^t_r(y)) = \partial_t d^s_{(\Omega_t)_r}(y)
$$
for $y \in \hat{\partial}( \Omega_t)_r$.

Thus we have showed that for a.e. $t \in [0, T]$, 
every $r \in (0, r_0]$ the equality holds
\begin{equation}
\label{aeTcomprmold}
\int_{\bR^2}w \varphi \,dz = 
\int_{\hat{\partial} (\Omega_t)_r }
(\varphi\circ P^t_r) \,  \partial_t d^s_{(\Omega_t)_r}(1- r\kappa_r^t) 
d{\cal H}^1
+ \int_{\Omega_t}aD \varphi Du \,dz.
\end{equation}

In the next lemma we show that in the first integral at the
right-hand side of (\ref{aeTcomprmold}) we can integrate
over the whole boundary $\partial (\Omega_t)_r$.

\begin{lemma}
\label{Int0}
Let $r \in (0, r_0]$ be such that properties (\ref{assm1}), 
(\ref{assm2}), (\ref{assm3}) are satisfied. Then for
 a.e. $t \in (0, T)$ 
$$
\int_{\partial (\Omega_t)_r \, \cap \{\gamma_{(\Omega_t)_r} \leq r \} } 
\varphi \partial_t d^s_{(\Omega_t)_r}(1- r\kappa_r^t) d{\cal H}^1 =
0
$$
\end{lemma}
\begin{proof}
We will show that for a.e. $t$
\begin{equation}
\label{intgrnd0}
\partial_t d^s_{(\Omega_t)_r}(y)(1- r\kappa_r^t(y)) = 0
\;\;\; {\cal H}^1 \;\mbox{ a.e. on }
\partial (\Omega_t)_r \, \cap \{\gamma_{(\Omega_t)_r} \leq r \}.
\end{equation}

By (\ref{assm2}) for a.e. $t$ the surface ${\bf \Gamma}_r$
  is (2,1)
differentiable ${\cal H}^1$ a.e. on $\partial (\Omega_t)_r$.
Fix such $t$. Then in (\ref{intgrnd0}) we can consider
only $y$ at which ${\bf \Gamma}_r$ is (2,1) differentiable.
Fix such $y$.

By Proposition \ref{Proposition 2}, 
$$
\gamma_{(\Omega_t)_r} \geq r \;\; \mbox{ on } \; \partial (\Omega_t)_r.
$$
Thus
$$
\gamma_{(\Omega_t)_r}(y) = r.
$$
It follows that
$$
\kappa_r^t(y) \leq \frac{1}{r}.
$$
Let $y' = P^t_r(y)$, i.e., $y'$ is the  point of 
$ \partial \Omega_t$ nearest to $y$. 
Consider 3 cases.

{\bf Case 1.} {\it $y$ is the 
unique point of 
$\partial (\Omega_t)_r$ for which $y'$ the nearest point on
$\partial \Omega_t$.} 

Then the calculations of \cite{EG}, Proposition 7.1,
Steps 2,3 imply that
$$
\kappa_r^t(y) = \frac{1}{r}.
$$
Thus we have (\ref{intgrnd0}) in this case.

{\bf Case 2.}{\it There exists $z \in \partial (\Omega_t)_r$, $z \neq y$,
such that $y' = P^t_r(z)$, and }
$$
B_r(y) \cap B_r(z) \neq \emptyset.
$$

Note that
$$
|y-y'| = |z-y'| = r.
$$
Then, denoting $w$ the point $\frac{y+z}{2}$,
we get
$$
w \in B_r(y) \cup B_r(z).
$$
We have 
$$
B_r(y),\; B_r(z) \subset \bR^2 \setminus \overline{\Omega}_t.
$$
Thus
$
w \in \bR^2 \setminus \overline{\Omega}_t,
$
and $y'$ is the  point of $\partial \Omega_t$ nearest to $w$.
Denote
$$
v=y' + \frac{r}{|w-y'|}(w-y').
$$
Then, by Condition \ref{Cond.1}, we get the following:
$$
v \notin \Omega_t, \;\;\;
y' \; is \; the\;unique \; point \; of\; \partial\Omega_t \; 
nearest \; to \; v. 
$$
But 
$$
|v-y'|=r.
$$
Thus we proved that
$$
v \in \partial (\Omega_t)_r.
$$
The points $y$ and $z$ divide
the circle $\partial B_r(y')$
 on two arcs, and the point $v$ is the middle point
of one of these arcs. Denote this arc $C_1$. Repeating
the above argument inductively and using continuity of
distance function, we prove that
$$
C_1 \subset \partial (\Omega_t)_r.
$$
But then 
$$
\kappa_r^t(y) = \frac{1}{r}.
$$
Thus we have (\ref{intgrnd0}) in the Case  2.

{\bf Case 3.}{\it There exists $z \in \partial (\Omega_t)_r$, $z \neq y$,
such that $y' = P^t_r(z)$, and }
$$
B_r(y) \cap B_r(z) = \emptyset.
$$

Then 
$$
\partial B_r(y) \cap \partial B_r(z) = {y'}.
$$
Since $\{ \Omega_t \}$ is an expanding family of sets, we have
\begin{equation}
\label{expand}
B_r(y),\; B_r(z) \subset \bR^2 \setminus \overline{\Omega}_{\tau}
\;\;\; \mbox{ for } \; \tau < t.
\end{equation}
If
\begin{equation}
\label{subCond1} 
y' \in \partial \Omega_{\tau^*} \;\;\;
\mbox{ for some } \;\; \tau^* < t,
\end{equation}
then the same is true for all $\tau \in [\tau^*, t]$, and
(\ref{expand}) implies that 
$$
y \in \partial (\Omega_{\tau})_r \;\;\; \mbox{ for all } \;\;
\tau \in [\tau^*, t].
$$
Then
$$
\partial_t d^{s}_{(\Omega_t)_{r}} (y) = 0.
$$
Thus we have (\ref{intgrnd0}) if (\ref{subCond1}) is satisfied.

The remaining case is 
\begin{equation}
\label{subCond2} 
y' \notin \partial \Omega_{\tau} \;\;\;
\mbox{ for all } \;\; \tau < t.
\end{equation}
Introduce a coordinate system $(x_1, x_2)$ on $R^2$ in which
$$
y'=(0,0),\;\; y=(0, r).
$$
Then
$$
z=(0, -r).
$$
Since $\{\Omega_t\}$ is a continuous  expanding family of sets,
the function
$$
\phi (\tau) = \dist(y', \Omega_{\tau}) 
$$
is continuous and nonincreasing, and
$$
\phi (t) = 0.
$$
Let $w_{\tau} \in {\cal N}_{\partial \Omega_{\tau}}(y')$. 
Then by (\ref{expand})
\begin{equation}
\label{wTau}
w_{\tau} \in \partial B_{\phi (\tau)}(0,0) \setminus
\left [ B_r(0,r) \cup B_r(0, -r)) \right ].
\end{equation}
In particular,
\begin{equation}
\label{convWtau}
w_{\tau} \rightarrow y'\;\;\; \mbox{as}\;\;\;  \tau \rightarrow t.
\end{equation}
By Condition  \ref{Cond.1} we get
\begin{equation}
\label{vBall}
B_r(v_{\tau}) \subset \bR^2 \setminus \overline{\Omega}_{\tau} \;\;\;
\mbox{ where } \;\; 
v_{\tau}= w_{\tau} + \frac{r}{| y' - w_{\tau}|}(y' - w_{\tau}).
\end{equation}
By (\ref{wTau}), (\ref{convWtau})  we see
that 
there exists a 
sequence 
$\tau_j \rightarrow t$ 
such that
$$
v_{\tau_j} \rightarrow  p, \;\; \mbox{where} \;\; p \;\; 
\mbox{is either} \;\;
(-r, 0) \;\; \mbox{or} \;\; (r, 0).
$$
Then by (\ref{vBall}) and continuity of the family
$\{ \Omega_t \}$ we conclude
$$
\mbox{either } \;\;B_{r} (-r, 0) 
\subset \bR^2 \setminus \overline{\Omega}_t
\;\;\mbox{ or }\;\;
B_{r}(r, 0) \subset \bR^2 \setminus \overline{\Omega}_t.
$$
Let in fact
$$
B_{r}(r, 0) \subset \bR^2 \setminus \overline{\Omega}_t.
$$
Then $y'=(0,0)$ is 
the  point
of $\partial \Omega_t$ nearest to the point $(r, 0)$. 
Thus
$$
(r, 0) \in \partial (\Omega_t)_r.
$$
But $y=(0, r)$, and thus
$$
B_r(y) \cap B_r(r, 0) \neq \emptyset.
$$
Thus  the points
$y$ and $(r, 0)$ satisfy the conditions of Case 2. Thus we get
$$
\kappa_r^t(y) = \frac{1}{r}.
$$
Case 3 is proved. Then Lemma \ref{Int0} is proved.
\end{proof}

The properties (\ref{assm1}), (\ref{assm2}), (\ref{assm3}) are satisfied
for a.e.  $r \in (0, r_0]$. Then there exists a decreasing
sequence $r_i \rightarrow 0 $, where $i=1, 2, ...$ and
 $r_i \in (0, r_0]$, such that (\ref{assm1}), 
(\ref{assm2}), (\ref{assm3}) are satisfied for each $r_i$.
Note that by (\ref{assm3}), the function
$$
(y,t) \rightarrow \kappa_{r_i}^t(y)
$$
defined  on ${\bf \Gamma}_{r_i}$ is ${\cal H}^2$ measurable.

Now we can integrate (\ref{aeTcomprmold}) by $t$
 and use Lemma  \ref{Int0} to get
\begin{eqnarray}
\nonumber
\int_0^T \int_{\bR^2}w \varphi \,dzdt 
 & = &  \int_0^T \int_{\partial (\Omega_t)_{r_i} }
(\varphi\circ P^t_{r_i}) \,  \partial_t d^s_{(\Omega_t)_{r_i}}
(1-r_i \kappa_{r_i}^t) 
d{\cal H}^1\,dt 
\\
&& +  \int_0^T \int_{\Omega_t}aD \varphi Du \,dzdt.
\label{0Tcomprmold}
\end{eqnarray}

It follows that the first integral at the right-hand side does
not depend on $r_i$. Thus it is enough to compute
the limit as $i \rightarrow \infty$.
We will prove the following:
\begin{equation}
\label{IndFunctEqual}
\lim_{i\rightarrow\infty}
\int_0^T \int_{\partial (\Omega_t)_{r_i} }
(\varphi\circ P^t_{r_i}) \,  \partial_t d^s_{(\Omega_t)_{r_i}}
(1- r_i \kappa_{r_i}^t) 
d{\cal H}^1\,dt =
-\int_0^T \int_{\bR^2} \chi_{\Omega_t} \partial_t\varphi\,dxdt
\end{equation}

We first prove that such equality is true if
the boundary satisfies additional regularity assumptions.
\begin{lemma}
Let $\{ \Omega_t \}$ satisfy (\ref{assm1}) - (\ref{assm3})
with $r=0$.
Let $V(x,t)$ be the
outer normal velocity of $\Gamma_t$,
defined by (\ref{outNormVelDef}) at every point $(x,t)$ of
differentiability of ${\bf \Gamma}$ and defined as 0 at
all points where ${\bf \Gamma}$ is not differentiable. Then
for every $\varphi \in C^1_c(\bR^2 \times (0,T))$
\begin{equation}
\label{kinetimatic}
\int_0^T \int_{\partial \Omega_t }
\varphi V d{\cal H}^{1} dt =
-\int_0^T \int_{\bR^2} \chi_{\Omega_t} \partial_t\varphi\,dxdt
\end{equation}
\end{lemma}
\begin{proof}
Let the set
$E$ be defined by (\ref{defE}). Then 
by (\ref{assm1})
$E$ has locally finite
perimeter. By (\ref {assm2})
$$
{\cal H}^2 (\partial E \setminus \partial_* E )= 0,
$$
 where 
$\partial_* E$ is the reduced boundary of $E$. Let 
$\Phi: \bR^2 \times \bR^1 \rightarrow \bR^2 \times \bR^1$
be defined by
$$
\Phi(x,t) = (0,...,0,\varphi(x,t))
$$
Let $\nu_E(x,t)$ be the measure-theoretical outer normal to
${\bf \Gamma}$ at $(x,t) \in {\bf \Gamma}$. Then
by Green-Gauss theorem for sets with finite perimeter
(\cite{EGar}, section 5.8)
$$
\int_0^T \int_{\bR^2} \chi_{\Omega_t} \partial_t\varphi\,dxdt =
\int_E \divg \Phi\,dxdt = \int_{\bf \Gamma} \Phi\, \nu_E\, d{\cal H}^2,
$$
where 
$$
\divg \Phi =\sum_{i=1}^{3} \frac{\partial \Phi_i}{\partial x_i},\;\;\;
\mbox{ where }\;x_{3}=t.
$$
At every point $(x,t)$ of differentiability of ${\bf \Gamma}$
 we have 
$$
\nu_E(x,t) = - D_{(x,t)}d^s_E(x,t),
$$
where $d^s_E$ is the signed distance to the boundary of the set
$E$ in the $(x,t)$-space.  
Thus we get
$$
\int_E \partial_t\varphi\,dxdt = 
- \int_{\bf \Gamma}\varphi\, \partial_t d^s_E \,d{\cal H}^2.
$$
 Let $f: {\bf \Gamma} \rightarrow \bR^1$ be the mapping
defined by
$$
f(x,t) = t
$$
Then $f$ is Lipschitz, and $f$ and ${\bf \Gamma}$ 
are differentiable at ${\cal H}^2$ a.e. point 
$(x,t) \in {\bf \Gamma}$. At such point $(x,t)$ the gradient of $f$ is
a linear mapping $Df(x,t):\,T_{x,t}{\bf \Gamma} \rightarrow \bR^1$,
where $T_{x,t}{\bf \Gamma}$ is the tangent to ${\bf \Gamma}$
at $(x,t)$ space.
Let $e_1, e_2$
be such orthonormal basis in $\bR^2$ that $e_2$ is the 
inner normal to $\partial \Omega_t$ at $x$.
Let $\tau$ be the unit vector in the $t$-direction. Then the vectors 
$$
e_1,\; 
\tilde{e}_2=-\partial_t d^s_E(x,t) e_2 + |D_x d^s_E(x,t)|\tau
$$
form an orthonormal basis in $T_{x,t}{\bf \Gamma}$.
We calculate:
\begin{eqnarray}
Df(x,t)e_1 &=& 0,
\nonumber \\
Df(x,t)\tilde{e}_2 &=& |D_x d^s_E(x,t)|.
\nonumber
\end{eqnarray}
Thus
$$
|Df(x,t)|= |D_x d^s_E(x,t)|.
$$
Also, the following relation holds:
$$
D_{x,t} d^s_{\Omega_t}(x,t) = \frac{1}{|D_x d^s_E(x,t)|}D_{x,t} d^s_E(x,t).
$$
Now, applying formula 3.2.22 of \cite{Federer}
(which is applicable by (\ref{assm1})), we get:
\begin{eqnarray}
\int_{\bf \Gamma}\varphi\, \partial_t d^s_E \,d{\cal H}^2
& = & \int_{\bf \Gamma}\varphi \, \partial_t d^s_{\Omega_t} 
|Df(x,t)| \,d{\cal H}^2
\nonumber \\
&= & \int_0^T \int_{f^{-1}({t})}\varphi \, \partial_t d^s_{\Omega_t}
\,d{\cal H}^{1}dt
\nonumber \\
&= & \int_0^T \int_{\partial \Omega_t} \varphi \, V
\,d{\cal H}^{1}dt.
\nonumber
\end{eqnarray}
The lemma is proved.
\end{proof}

Now we can prove (\ref{IndFunctEqual}). 
Each $\{ (\Omega_t)_{r_i} \}$ satisfies (\ref{assm1}) - (\ref{assm3}).
Thus we have
\begin{equation}
\label{approx}
\int_0^T \int_{\partial (\Omega_t)_{r_i} }
\varphi V^{(r_i)} d{\cal H}^{1} dt =
-\int_0^T \int_{\bR^2} \chi_{(\Omega_t)_{r_i}} \partial_t\varphi\,dxdt,
\end{equation}
where $V^{r_i}$ is velocity of $(\Omega_t)_{r_i}$.
It follows from Proposition \ref{Proposition 4.1} that
${\cal L}^3({\bf \Gamma}) = 0$, and thus
$$
\chi_{(\Omega_t)_{r_i}}(x) \rightarrow \chi_{\Omega_t}(x)
\;\;\;\; \mbox{for a.e. } (x,t) \in \bR^n\times(0,T).
$$
Then by Dominated Convergence Theorem the right-hand side of
(\ref{approx}) converges to the right-hand side of (\ref{IndFunctEqual})
as $i \rightarrow \infty$.
Thus it remains to prove that the left-hand side of
(\ref{approx}) and the left-hand side of (\ref{IndFunctEqual})
converge to the same limit.

Let $P^t_{r, r_i }: \partial (\Omega_t)_{r} \rightarrow 
 \partial (\Omega_t)_{r_i}$
be the nearest point projection (well-defined by
Proposition \ref{Proposition 2}). $P^t_{r, r_i }$ is
a Lipschitz map by Propositions \ref{Proposition 2}, 
\ref{Proposition 3}.
Then  using Lemma \ref{jacobian lemma}
and identity (\ref{relNbhd}), the difference between 
the  left-hand sides of (\ref{approx}) and (\ref{IndFunctEqual})
can be transformed to
\begin{eqnarray}
\label{lhs1}
&& \int_0^T \int_{\partial (\Omega_t)_{r} }
(\varphi\circ P^t_{r, r_i} - \varphi\circ P^t_{r}) \,  
V^{r_i}\circ P^t_{r, r_i} \,
(1- r \kappa_{r}^t) 
d{\cal H}^1\,dt + \\
&& r_i\int_0^T \int_{\partial (\Omega_t)_{r} }
\varphi\circ P^t_{r, r_i} \,  
V^{r_i}\circ P^t_{r, r_i}\, \kappa_{r}^t\,
d{\cal H}^1\,dt = I_{1, i} + r_i I_{2, i}.
\nonumber
\end{eqnarray}

We have $|V^{r_i}|<C$, where $C$ does not
depend on $i$, and $|\kappa_{r}^t| \leq \frac{1}{r}$. 
Using the fact that
$P^t_{r, r_i} \rightarrow P^t_r$ as $r_i \rightarrow 0$
 and Dominated 
Convergence Theorem we see that $I_{1, i} \rightarrow 0$. 
We also have $|I_{2, i}| < C$.
Thus the expression (\ref{lhs1})
converges to zero. Thus (\ref{IndFunctEqual}) is proved.

The equalities (\ref{0Tcomprmold}) and (\ref{IndFunctEqual})
imply (\ref{balEqComprMold}).

{\it Step 3}.

The family $\{\Omega_t\}$ is expanding. Thus the left-hand side of
(\ref{IndFunctEqual}) defines a nonnegative linear functional
of $\varphi \in C^{\infty}(\bR^n \times [0,T])$. 
By \cite{EGar}, Chapter 1.8, Corollary 1, it follows 
from (\ref{IndFunctEqual}) that
$\partial_t \chi_{\Omega_t}$ is a nonnegative
Radon measure. Thus the mass balance equation (\ref{balEqComprMold})
can be rewritten as 
$$
\int_0^T \int_{\bR^2} w \varphi\, dxtdt + 
\int_{\bR^2 \times (0,T)}\varphi \,dw_t
-\int_0^T \int_{\bR^2} aDuD\varphi\, dxtdt = 0.
$$
From this equation it follows that for a.e. $t \in [0,T]$,
every $\varphi \in C^1_c(\bR^2)$
$$
\int_{\bR^2} w(\cdot, t) \varphi\, dx + 
\int_{\bR^2}\varphi \,dw_t(\cdot, t)
-\int_{\bR^2} a(\cdot, t)Du(\cdot, t)D\varphi\, dx = 0.
$$
Now, repeating argument of Section \ref{Th.3}, we conclude the
proof of Theorem \ref{Theorem 1-CM}.
\end{proof}

\appendix

\section{Appendix}
\label{Apndx 1}

Let $\Omega$ be an open set.
We write $d(x)$
for $d_{\Omega}(x)$ below. The purpose of this section is to
prove that the gradient of $d_{\Omega}$  is locally Lipschitz at
$x \in \Omega \setminus {\cal R}$ and give an estimate of the
Lipschitz constant in the terms of the distance between $x$ and
endpoints of the distance ray $R_x$. 

Proposition \ref{Proposition 3Ap} should
be compared with Proposition 4.1 of \cite{EG} and
with inequality 4.8(8) of \cite{FedCMeas}. 

In the inequality (\ref{1.1a})
below the only assumption regarding the point $x_1$ is that
it is close enough to $x$. In particular it is possible that 
$x_1 \in {\cal R}$. In
the inequality 4.8(8) of \cite{FedCMeas} the conditions
on $x$ and $x_1$ are symmetric and exclude the possibility that
$x_1 \in {\cal R}$.

Proposition \ref{Proposition 3Ap} 
improves the estimates of Proposition 4.1 of \cite{EG} 
in the following.
Two quantities are estimated explicitly
in Proposition \ref{Proposition 3Ap}:
the local Lipschitz constant of $Dd(\cdot)$ at $x$ and the size 
of the neighborhood of $x$ in which the estimates (\ref{1.1a})
and (\ref{1.1}) hold. 

\begin{proposition}
\label{Proposition 3Ap}
There exist constants $C$ and $M$ depending only
on $ n$ such that the following is true.
Let $\Omega \in R^n$ be an open set. 
Let $\varepsilon > 0$,  $x\in \Omega 
\backslash {\cal R}$, and let 
\begin{equation}
\label{ptsAssm}
d_{\Omega}(x) \geq M\varepsilon, \;\;\;
\gamma_{\Omega}(x) - d_{\Omega}(x) \geq M\varepsilon.
\end{equation}
Then by Remark \ref{remUnNr}
there exists a unique ray $R_x$. Denote $y$ and $v$ the lower and 
upper ends of $R_x$, i.e., $y=R_x \cap \partial \Omega$ and
$v=R_x \cap {\cal R}$.
Then for every
 $x_1\in \Omega$ satisfying $\mid x-x_1 \mid <\varepsilon$,  
the inequality holds
\begin{equation} 
\label{1.1a}
\mid y-y_1\mid \leq C\left(1+\frac{ \mid x-y \mid}{ \mid x-v \mid}\right) 
\mid x-x_1 \mid,
\end{equation}
where $y_1 \in {\cal N}_{\partial\Omega}(x_1)$.
If in addition the function $d(\cdot)$ is differentiable at $x_1$ then 
\begin{equation} 
\label{1.1}
|Dd(x)-
Dd(x_1)|\leq \frac{C }{ \varepsilon }|x-x_1|.
\end{equation}
\end{proposition}
\begin{proof} 
We first prove the inequalities (\ref {1.1a}), (\ref{1.1}) assuming
that
\begin{equation}
\label{convCond}
\mid x- v \mid \geq 
\mid x- y \mid.
\end{equation}

Denote $d:=d_{\Omega}(x)$.
By  (\ref{convCond}), we can find a point $O$ on the interval of $R_x$
between $v$ and $x$ such that $\mid x-O\mid =  d$. Let 
$y \in {\cal N}_{\partial \Omega}(x)$, then 
we also have $y \in {\cal N}_{\partial \Omega}(O)$.
In the calculations below $C$ will denote
different constants depending only on $n$. We assume that 
\begin{equation}
\label{Mgt10}
M>10.
\end{equation}

Choose $x_1 \in B_{\varepsilon}(x)$. Let 
$\hat{x}_1$ be projection of $x_1$ onto $R_x$,
then 
$$|x-\hat{x}_1|\leq \varepsilon \leq {1\over 10}d$$
by (\ref{ptsAssm}), (\ref{Mgt10}).
Thus 
$\hat{x}_1$ lies
between $y$ and $O$ on $R_x$. Let $d_1=|\hat{x}_1-O|$, then 
\begin{equation}
\label{1.2}
{11\over 10}\geq {d\over {d_1}}\geq{9\over 10},\; \; \; \;|d-d_1|
<\varepsilon
\end{equation}
Introduce the following coordinate system in $R^n$: let the point $O$
be the origin, let $e_n=\frac{x-O}{|x-O|}$ (thus $e_n$ is the unit vector 
along the ray $R_x$), let $e_1,...,e_{n-1}, e_n$ be an orthonormal basis in
$\bR^n$.  Then in these coordinates
$$x=(0,d)_,\;\;\;\hat{x}_1=(0,d_1),\;\;\;y=(0,2d),$$
where $0 \in \bR^{n-1}$.  We also have  $x_1=(x'_1,d_1),$
where $x'_1 \in \bR^{n-1}$. Let $y_1 \in \partial
\Omega$ be such that $|x_1-y_1|=d(x_1)$. Let the coordinates of
$y_1$ be $(y'_1,y_{1,n})$ where $y'_1 \in \bR^{n-1}$, $ y_{1,n} \in
\bR^1$.  Then
\begin{equation}
\label{1.3}
|x_1-y_1|\leq|x_1-y|
\end{equation}
or
\begin{equation}
\label{1.4}
|x'_1-y'_1|^2+|d_1-y_{1,n} |^2\leq|x'_1|^2+|d_1-2d|^2
\end{equation}

\begin{claim}
\label{Claim 1}
There exists $M_1$,
depending on $n$, such that for any $ \delta \in (0,2)$ 
the following is true:  
\begin{equation}
\label{claimCond}
\mbox{if} \;\;d>\frac{M_1}{\delta^2}\varepsilon\;\; 
\mbox{then} \;\;\;|y'_1|<\delta d.
\end{equation}
\end{claim}

\begin{proof}
Note that 
\begin{equation}
\label{1.41}
|x'_1| < \varepsilon, \;\;\;\; |y'_1|<2d.
\end{equation}
Indeed, the first inequality is true since $x_1 \in B_{\varepsilon}(x)$. 
To prove the second inequality of (\ref{1.41})
 we use (\ref{1.4}) and (\ref{1.2}) to get
$$
|y'_1- x'_1|^2 \leq \varepsilon^2+|d+\varepsilon|^2
<(d+\sqrt{2}\varepsilon)^2,
$$ 
so 
\begin{equation}
\label{1.5}
|y'_1|\leq |x'_1|+d+\sqrt{2}\varepsilon<d+3\varepsilon<2d
\end{equation}
since $d\geq 10\varepsilon$. Thus  (\ref{1.41}) is proved.

Suppose that the assertion of Claim \ref{Claim 1} is false,
that is  
\begin{equation}
\label{wrongIneq}
|y'_1|\geq \delta d.
\end{equation} 
By (\ref{claimCond})
\begin{equation}
\label{deltaD}
\delta d > {M_1 \over \delta}\varepsilon > 2\varepsilon
\end{equation}
if $M_1 > 4$.
Then we get from (\ref{1.4}), (\ref{1.41}),  (\ref{1.2}),
(\ref{wrongIneq})
\begin{eqnarray}
(y_{1,n} -d_1)^2 & \leq & 
\varepsilon^2+(2d-d_1)^2-(|y'_1|-|x'_1|)^2
\nonumber \\
 & \leq &
\varepsilon^2+({11 \over 10}d)^2-(\delta d-\varepsilon)^2
\label{1.6} \\
& < & 4d^2(1-
\frac{\delta^2}{4})-2d^2+2\delta d\varepsilon 
\nonumber \\
 & \leq &
 4d^2(1-
\frac{\delta^2}{4})
\nonumber
\end{eqnarray}
The last inequality follows from
 $-2d^2+2\delta d\varepsilon < 0$ which holds
by (\ref{claimCond}) if
$\delta < 2$ and $ M_1>8$.

Consider two cases: \\
\underline{Case 1}:$\;\; y_{1,n}  \geq d_1$.
Denote by $\alpha$ the angle between vectors $y_1 - x_1$ and 
$(x'_1,y_{1,n})-x_1$.
In the triangle $x_1, y_1,\;(x'_1,y_{1,n})$ the side $(x'_1,y_{1,n})-
y_1$ is orthogonal to $x_1-(x'_1,y_{1,n})$, since $(x'_1,y_{1,n})-
y_1=(x'_1-y'_1,0)$, and 
$x_1-(x'_1,y_{1,n})=(0,d_1-y_{1,n})$. Thus we get
$$
\tan\alpha=\frac{|y'_1-x'_1|}{y_{1,n}  - d_1}.
$$
But 
$$|y'_1-x'_1|\geq||y'_1|-|x'_1||=|\delta d-
\varepsilon|=\delta d-\varepsilon
$$
by (\ref{deltaD}).  By (\ref{1.6}) we get
$$
\tan\alpha \geq \frac{\delta  d-\varepsilon}{2d\sqrt{1-
\frac{\delta^2}{4}}}=\frac{\delta}{2\sqrt{1-
\frac{\delta^2}{4}}}(1-{1\over \delta}\cdot{\varepsilon \over d})\geq
\frac{\delta}{4\sqrt{1-\frac{\delta^2}{4}}}.
$$
where the last inequality follows from (\ref{deltaD}).
Thus,
using the condition $\delta \in (0, 2)$ we get 
\begin{equation}
\label{1.7}
|\cos\alpha|=\frac{1}{\sqrt{1+\tan^2\alpha}}\leq
\sqrt{\frac{4(4-\delta^2)}{16-
3\delta^2}} \leq 1-
{\delta^2\over {2(16 - 3\delta^2)}}.
\end{equation}

Consider the triangle $x_1,y_1,\;(x'_1,0)$.
The angle between the vectors $(x'_1,0) - x_1$ and $y_1-x_1$ is $\pi-
\alpha$, and so we get
\begin{equation}
\label{ThCos}
|y_1-(x'_1,0)|^2=|x_1-(x'_1,0)|^2+|y_1-x_1|^2-2|x_1-(x'_1,0)||y_1-
x_1|\cos(\pi-\alpha).
\end{equation}
Since $x_1=(x'_1, d_1)$, we have
$$
|x_1-(x'_1,0)|=d_1.
$$
From (\ref{1.4}) we get
$$
|y_1-x_1|\leq\sqrt{\varepsilon^2+(2d-d_1)^2}.
$$
Then by (\ref{ThCos}), (\ref{1.7})
we get
\begin{eqnarray}
|y_1-(x'_1,0)|^2
& \leq &
 d^2_1 +[\varepsilon^2+(2d-
d_1)^2]+2d_1\sqrt{\varepsilon^2+(2d-d_1)^2}\,(1-{\delta^2 \over 32})
\nonumber \\
& \leq &
 d^2_1+(2d-d_1)^2+2d_1(2d-d_1)(1+\frac{\varepsilon^2}{2(2d-
d_1)^2})(1-{\delta^2 \over 32})+\varepsilon^2
\nonumber \\
& \leq  &
4d^2+2d_1\frac{\varepsilon^2(1-{\delta^2\over 32})}{2(2d-
d_1)}+\varepsilon^2-2d_1(2d-d_1){\delta^2 \over 32}
\nonumber \\
& \leq & 
4d^2+d_1\frac{\varepsilon^2}{2d-d_1}+\varepsilon^2-{1\over
16}d_1(2d-d_1)\delta^2.
\nonumber
\end{eqnarray}
Thus we see that
\begin{equation}
\label{mainClaimEst}
|y_1-(x'_1,0)|^2 < (2d-\varepsilon)^2,
\end{equation}
provided
\begin{equation}
\label{lastInq}
d_1\frac{\varepsilon^2}{2d-d_1}-{1\over 16}d_1(2d-
d_1)\delta^2<-4d\varepsilon.
\end{equation}
Using  (\ref{1.2}) we see that (\ref{lastInq}) is satisfied if
\begin{equation}
\label{nearContr}
d^2\delta^2>C(\varepsilon^2+d\varepsilon)
\end{equation}
where $C$ is a large enough constant.  
In its turn, (\ref{nearContr}) is true if 
$d>{M_1\over \delta^2}\varepsilon$, 
where $M_1$ is large enough.
Thus, for such $d$, the inequality (\ref{mainClaimEst}) holds, and so
$$
|y_1-(x'_1,0)|<2d-\varepsilon.
$$
Then 
$$
|y_1|\leq|x'_1|+|y_1-(x'_1,0)|<2d=|y|.
$$
This contradicts the fact
that $y \in {\cal N}_{\partial \Omega}(O)$. 
Thus (\ref{wrongIneq}) is false.
Thus Case 1 of Claim \ref{Claim 1} is proved.

\underline{Case 2:} $y_{1,n}  \leq d_1.$

Using this condition, we get
\begin{eqnarray}
|(x'_1,0)-y_1|^2 & = & 
|(x'_1,0)-x_1+x_1-y_1|^2=|-(0,d_1)+(x'_1-
y'_1,\;d_1-y_{1,n})|^2
\nonumber \\
& = & d^2_1+|x_1-y_1|^2+2(y_{1,n}  
-d_1)d_1
\nonumber \\
& \leq & d^2_1+|x_1-y_1|^2\leq
d^2_1 +\varepsilon^2+(2d-d_1)^2
\nonumber \\
& = & 4d^2-2d_1(2d-d_1)+\varepsilon^2<(2d-\varepsilon)^2,
\nonumber \\
& & \mbox{ if }\ d_1(2d-d_1)>2d\varepsilon.
\nonumber 
\end{eqnarray}
The last inequality is in fact true by
(\ref{1.2}) and $d>10\varepsilon$.  Then, as in
the Case 1, we arrive to a contradiction. Thus
Case 2 cannot happen for $d$ and $ \varepsilon$ defined in 
Claim \ref{Claim 1}.  Thus,
Claim \ref{Claim 1} is proved.
\end{proof}

Now the rest of the proof of part (a) of 
Proposition \ref{Proposition 3} 
follows the
proof of  Proposition 4.1 in \cite{EG}.

The main step is to show that there exist  
constants $C$ and $ M$ depending only on $n$ such that the 
inequality $d>M\varepsilon$ implies
\begin{equation}
\label{1.11}
|y-y_1|\leq C|\hat{x}_1-x_1|.
\end{equation}

Let $\delta \in (0,1)$ be chosen later, and let $d>{M_1\over
\delta^2}\varepsilon$ where $M_1$ is the constant from
Claim \ref{Claim 1}. Then by Claim \ref{Claim 1} we get
\begin{equation}
\label{RightIneq}
|y'_1|<\delta d.
\end{equation}
Since $y\in {\cal N}_{\partial \Omega}(O)$, we get $|y_1|>|y|=2d$, 
and so  
$$
|y'_1|^2+|y_{1,n}  |^2\geq
4d^2.
$$
Thus we estimate using (\ref{RightIneq})
\begin{eqnarray}
y_{1,n} & \geq & \sqrt{4d^2-|y'_1|^2}  =  2d\sqrt{1-
\frac{|y'_1|^2}{4d^2}}\geq 2d(1-\frac{|y'_1|^2}{8d^2}-C{|y'_1|^4
\over {d^4}})
\nonumber \\
\label{1.899}
& \geq & 2d-\frac{|y'_1|^2}{4d}-C\frac{|y'_1|^4}{4d^3}\geq 2d-
\frac{|y'_1|^2}{4d}(1+C\delta^2),
\end{eqnarray}
where all constants $C$ do not depend on $y'_1$ and $d$.  So
$$
0<2d-d_1=y_{1,n}  -d_1+2d-y_{1,n}  \leq (y_{1,n}  -
d_1)+\frac{|y'_1|^2}{4d}(1+C\delta^2).
$$
Now we can estimate
\begin{eqnarray}
(2d-d_1)^2 & \leq & (y_{1,n}  -d_1)^2+(y_{1,n}  -
d_1)\frac{|y'_1|^2}{2d}(1+C\delta^2)+\frac{|y'_1|^4}{16{d^2}}
(1+C\delta^2)
\nonumber \\
\label{1.9}
& \leq & (y_{1,n}  -d_1)^2+(y_{1,n}  -
d_1)\frac{|y'_1|^2}{2d}(1+C\delta^2)+C\delta^2 |y'_1|^2
\end{eqnarray}
where we again used (\ref{RightIneq}).
Using (\ref{1.4}), (\ref{1.2}) and (\ref{1.41}) 
we get the following estimate
\begin{eqnarray}
|y_{1,n}  -d_1| & \leq & \sqrt{\varepsilon^2 +(2d-d_1)^2}=(2d-
d_1)\sqrt{1+\frac{\varepsilon^2}{(2d-d_1)^2}}
\nonumber \\
& \leq & (2d-
d_1)(1+\frac{\varepsilon^2}{2(2d-d_1)^2})
\leq(2d-d_1)(1+C\frac{\varepsilon^2}{d^2})
\nonumber \\
&\leq & (2d-
d_1)(1+C\delta^4)\;\;\;\mbox{since}\;\;d\geq
\frac{M_1}{\delta^2}\varepsilon. \nonumber
\end{eqnarray}
Thus, by (\ref{1.9}):
\begin{eqnarray}
(2d-d_1)^2 & \leq & (y_{1,n}  -d_1)^2+ 
(2d-d_1)\frac{|y'_1|^2}{2d}(1+C\delta^4)+C|y'_1|^2\delta^2
\nonumber \\
& \leq & (y_{1,n}  -
d_1)^2+\Theta|y'_1|^2,
\label{1.10}
\end{eqnarray}
where
$$
\Theta=(1-\frac{d_1}{2d})(1+C\delta^4)+C\delta^2.
$$
Using (\ref{1.2}) we estimate
$$
\Theta \leq{11\over
20}(1+C\delta^4)+C\delta^2<1,
$$
if $\delta>0$ is small enough.
Fix such $\delta$. Let $M_1$ be the constant from 
Claim \ref{Claim 1} defined by this
$\delta$.
Let $M=\frac{M_1}{\delta^2}$ be our choice of the constant $M$
in the inequalities (\ref{ptsAssm}).
Then
the  inequalities (\ref{ptsAssm})
imply (\ref{1.10}) with $\Theta<1$.  Thus (\ref{1.10}) 
and (\ref{1.4}) give:
\begin{equation}
\label{1.10.2}
|y'_1-x'_1|^2\leq |x'_1|^2+\Theta |y'_1|^2.
\end{equation}
Now we estimate:
$$
\mid y'_1 \mid \leq  \mid y'_1- x'_1 \mid + \mid x'_1 \mid
$$
and from this, using (\ref{1.10.2}), we get
\begin{eqnarray}
\mid y'_1 \mid^2 & \leq & \mid y'_1- x'_1 \mid^2 + 
2 \mid y'_1- x'_1 \mid  \mid x'_1 \mid + \mid x'_1 \mid^2 
\nonumber \\
& \leq & (1+\mu) \mid y'_1- x'_1 \mid^2 + (1+\frac{C}{\mu}) 
\mid x'_1 \mid^2
\nonumber \\
& \leq & (1+\mu)( |x'_1|^2+\Theta |y'_1|^2)+ 
(1+\frac{C}{\mu}) \mid x'_1 \mid^2.
\nonumber
\end{eqnarray}
Choose $\mu>0$ so small that $(1+\mu)\Theta<1$, then we get
\begin{equation}
\label{1.10.3}
\mid y'_1 \mid \leq C \mid x'_1 \mid=C\mid x_1 - \hat{x}_1\mid.
\end{equation}

It remains to estimate $\mid y_{1,n} - y_n \mid$ where $y_n=2d$.

From (\ref{1.899}) we get 
\begin{equation}
\label{1.10.4}
2d-y_{1,n}  \leq C \frac{\mid y'_1 \mid^2}{d} \leq C\mid y'_1 \mid,
\end{equation}
since $\mid y'_1 \mid < \delta d$.

Since $x_1=(x'_1, d_1)$, $d(x_1)=\mid x_1 - y_1 \mid$, and
$d(\hat{x}_1)) = \mid \hat{x}_1 -y \mid = 2d-d_1$, we calculate
(using the fact that $\mbox{Lip}[d(\cdot)]=1$)
$$
y_{1,n}  -d_1 \leq \mid y_1-x_1 \mid =
2d-d_1+ d(x_1) - d(\hat{x}_1) \leq
2d-d_1+\mid x_1- \hat{x}_1 \mid
$$
and so
\begin{equation}
\label{1.10.5}
y_{1,n}  - 2d \leq  \mid x_1- \hat{x}_1 \mid.
\end{equation}
The estimates (\ref{1.10.3}),  (\ref{1.10.4}), and (\ref{1.10.5}) imply
(\ref{1.11}).

Now since 
$$
Dd(x_1)=\frac{x_1-y_1}{\mid x_1-y_1\mid}, \;\;\;
Dd(\hat{x}_1)=\frac{\hat{x}_1-y}{\mid \hat{x}_1-y\mid},
$$
we get the following estimate
$$
|Dd(x_1)-Dd(\hat{x}_1)|\leq \frac{|x_1-\hat{x}_1|+|y_1-
y|}{|\hat{x}_1-y|},
$$
then using the equality 
$|\hat{x}_1-y|=2d-d_1$ and the inequalities (\ref{1.2}) and
(\ref{1.10})
we deduce that
$$
|Dd(x_1)-Dd(\hat{x}_1)|\leq {C\over d}|x_1-\hat{x}_1|.
$$
Note that  $Dd(\hat{x}_1)=Dd(x)$ since the 
points $x$ and $\hat{x}_1$ are on the same
ray $R^1_x$. Now from $|\hat{x}_1-x_1|\leq |x-x_1|$, we 
get  (\ref{1.1}) in the case of  (\ref{convCond}).
Note also that  (\ref{1.11}) and 
inequality $|\hat{x}_1-x_1|\leq |x-x_1|$
imply (\ref{1.1a}).

In the case  $\mid x- v \mid < 
\mid x- y \mid $ we can  make the following reduction 
to the case (\ref{convCond}).
Let 
\begin{equation}
\label{defQ}
Q=\mid x- y \mid - 
 \mid x- v \mid,
\end{equation}
and let $\tilde \Omega$ be the set 
$ \{z  \in \Omega \;\; \mid \;\;  \dist (z, \partial \Omega) 
> Q  \}$. Then $B_{\varepsilon}(x) \subset 
\tilde{\Omega}$ since the function $\dist(\cdot, \partial \Omega)$ 
is Lipschitz with constant 1. For $z \in \tilde{\Omega} $ we have
\begin{equation}
\label{ReductDist}
\dist(z, \partial {\Omega)}) = \dist(z, \partial \tilde{\Omega}) + Q. 
\end{equation}
To see this consider
$z_1 \in {\cal N}_{\partial \Omega}(z)$, and let $z_2$ be
the point of intersection of the interval connecting $z$ and $z_1$ with
$ \partial \tilde{\Omega} $. Then $z_1 \in {\cal N}_{\partial \Omega}(z_2)$
 since $z_2$ lies in the relative interior of the distance ray $R_{z_1}$ 
that intersects $\partial \Omega$ in $z_1$.
Thus $\mid z_1 - z_2 \mid = Q$.
Now suppose that there exists a point $\hat{z}_2 \in \partial 
\tilde{\Omega}$ such that
$\mid z-\hat{z}_2\mid < \mid z-z_2\mid $. Let  
$\hat{z}_1 \in {\cal N}_{\partial \Omega}(\hat{z}_2)$, 
then $\mid \hat{z}_1 - \hat{z}_2 \mid = Q$.
 Then 
$$\mid z- \hat{z}_1 \mid \leq
\mid z - \hat{z}_2 \mid + \mid \hat{z}_2 - \hat{z}_1\mid < \mid z-z_2 \mid +
Q =
\mid z - z_1 \mid,
$$
a contradiction to the fact that 
$z_1 \in {\cal N}_{\partial \Omega}(z)$. This proves (\ref{ReductDist}). 
So the distance rays and the ridge set of the set $ \tilde{\Omega}$ are 
intersections of those of the set 
$\Omega$ with the set $ \tilde{\Omega}$.

Then  in (\ref{1.1}) we can consider $d(\cdot)$ as the 
distance to $\partial \tilde{\Omega}$,
so we can consider $\tilde{\Omega}$ instead of $\Omega$, and then the 
inequality (\ref{convCond}) is  satisfied.

To show (\ref{1.1a}), denote by 
$w$ and  $w_1$ the points of intersection of  
$\partial \tilde{\Omega}$ with the
intervals $(x,y)$ and $(x_1, y_1)$ respectively.  Then it follows from
(\ref{ReductDist})
$w \in {\cal N}_{\partial \tilde{\Omega}}(x)$ 
and $w_1\in {\cal N}_{\partial \tilde{\Omega}}(x_1)$, and 
$$
\mid w-y \mid = \mid w_1 - y_1 \mid =Q.
$$
Then it follows from (\ref{defQ}) that 
$$
\mid x-v \mid = \mid x-w\mid.
$$
Using (\ref{defQ})
and letting 
$
A= \frac{ \mid x-y \mid}{ \mid x-v \mid},
$
we get
\begin{eqnarray}
y &=& x+ A(w-x), \;\;\;
\nonumber \\ 
y_1 &=& x_1 + \left[1+(A-1)\frac{\mid x-w \mid}{\mid x_1-w_1\mid} \right]
( w_1 - x_1).
\nonumber
\end{eqnarray}
The inequality   (\ref{1.1a}) can be applied to the set $ \tilde{\Omega}$, so
$$
\mid w-w_1\mid \leq C \mid x-x_1\mid.
$$
Now we estimate:
\begin{eqnarray}
\mid y-y_1\mid & = & \left| x-x_1 + 
A\left[(w-x) - \frac{\mid x-w \mid}{\mid x_1-w_1\mid}(w_1 - x_1)\right] \right.
\nonumber \\
& & + \left.
\left(\frac{\mid x-w \mid}{\mid x_1-w_1\mid}-1\right) ( w_1 - x_1) \right| 
\nonumber \\
&\leq &  C(1+A)(\mid x-x_1\mid + \mid w-w_1\mid) 
\nonumber \\
&\leq &
 C(1+A) \mid x-x_1\mid.
\nonumber
\end{eqnarray}
This implies (\ref{1.1a}).
\end{proof}  

\section*{Acknowledgement} It is the author's pleasure to thank L. C. Evans
and R. F. Gariepy for their interest, insight and and encouragement.


\begin{thebibliography}{99}
\parskip=0pt
\small
\itemsep=0pt

\bibitem[Ar95]{ArCompMol} G. Aronsson. {\em
Asymptotic solutions of a compression molding problem}, 
Preprint LiTH-MATH-R-95-01, 1995,
Department of Mathematics, Linkoping University, 
Linkoping , Sweden.

\bibitem[AE]{AEComprMot} G. Aronsson, L. C. Evans. {\em
An asymptotic model for compression molding}, Forthcoming.

\bibitem[AEW96]{AEW} G. Aronsson, L. C. Evans, Y. Wu.
{\em Fast/Slow diffusion and growing sandpiles},
 J.Diff.Equat., 131 (1996), no. 2, 304-335.

\bibitem[BDM89]{BDM} T. Bhattacharya, E. DiBenedetto, J. Manfredi.
{\em Limits as $p\rightarrow \infty$ of $\Delta_pu=f$
and related extremal problems},
Rend. Sem. Mat. Univ. Politee, Torino, 1989.

\bibitem[EFG97]{EFG} L. C. Evans, M. Feldman, R. F. Gariepy.
{\em Fast/Slow diffusion and collapsing sandpiles},
J.Diff.Equat., 137(1997), 166-209


\bibitem[EGan]{EG}  L. C. Evans, W. Gangbo.
  { \em Differential equations methods for the Monge-Kantorovich 
mass transfer problem}, Preprint.


\bibitem[EGar92]{EGar}  L. C. Evans, R. F. Gariepy.
{ \em Measure theory and fine properties of functions}, CRC Press, 1992.

\bibitem[EH87]{EH}  W. D. Evans, D. J. Harris.
  { \em  Sobolev embeddings for generalized ridged domains}, 
Proc.London Math.Soc., 54(1987), 141-175.

\bibitem[Fed59]{FedCMeas}  H. Federer. {\em Curvature measures},
Trans. Am. Math. Soc., 93(1959), 418-491.

\bibitem[Fed69]{Federer} H. Federer. {\em Geometric measure theory}, 
Springer, 1969.
 
\bibitem[Fel]{Feld} M. Feldman. {\em Convex viscosity solutions of nonlocal
geometric motion of planar curves}, in preparation.

\bibitem[GT83]{GilTr} D. Gilbarg, N.Trudinger.
  { \em Elliptic Partial Differential Equations of second order (2nd Ed.)},
Springer, 1983.

\bibitem[Jan93]{Jan} U.Janfalk.
  { \em On certain problems concerning the $p$-Laplace operator},
Linkoping Studies in Science and Technology,
Dissertation \#326, Linkoping University, Sweden, 1993.

\bibitem[Kry87]{Krylov}  N. Krylov. 
{\em Nonlinear elliptic and parabolic equations
of second order}, Reidel, Dordrecht, 1987.


\bibitem[Son93]{S}  H. M. Soner. {\em Motion of a set by a curvature of its
boundary},
J.Diff.Equat., 101(1993), 313-372.


\end{thebibliography}
\end{document}